\documentclass[%
 aip,
 amsmath,amssymb,
preprint,%
]{revtex4-1}

\usepackage{graphicx}
\usepackage{dcolumn}
\usepackage{bm}
\usepackage{xcolor}

\usepackage[utf8]{inputenc}
\usepackage[T1]{fontenc}
\usepackage{mathptmx}

\begin{document}
\frenchspacing
\preprint{}
 
\title[]{Systematic comparison between the generalized Lorenz equations and DNS in the two-dimensional Rayleigh-B\'enard convection}

\author{Junho Park}
\email{ad5486@coventry.ac.uk}
\affiliation{Fluid and Complex Systems Research Centre, Coventry University, Coventry CV1 5FB, UK
}%
\author{Sungju Moon}%
\affiliation{School of Earth and Environmental Sciences, Seoul National University, Seoul 08826, South Korea
}%
\author{Jaemyeong Mango Seo}
\affiliation{%
Max Planck Institute for Meteorology, Bundesstra\ss e 53, 20146 Hamburg, Germany
}%
\author{Jong-Jin Baik}
\affiliation{School of Earth and Environmental Sciences, Seoul National University, Seoul 08826, South Korea
}%

\date{\today}

\begin{abstract}
The classic Lorenz equations were originally derived from the two-dimensional Rayleigh-B\'enard convection system considering an idealised case with the lowest order of harmonics. 
Although the low-order Lorenz equations have traditionally served as a minimal model for chaotic and intermittent atmospheric motions, even the dynamics of the two-dimensional Rayleigh-B\'enard convection system is not fully represented by the Lorenz equations, and such differences have yet to be clearly identified in a systematic manner. 
In this paper, the convection problem is revisited through an investigation of various dynamical behaviors exhibited by a two-dimensional direct numerical simulation (DNS) and {the generalized expansion of the Lorenz equations (GELE)} derived by considering additional higher-order harmonics in the spectral expansions of periodic solutions. 
Notably, the {GELE} allows us to understand how nonlinear interactions among high-order modes alter the dynamical features of the Lorenz equations including fixed points, chaotic attractors, and periodic solutions. 
It is verified that numerical solutions of the DNS can be recovered from the solutions of {GELE} when we consider the system with sufficiently high-order harmonics. 
At the lowest order, the classic Lorenz equations are recovered from  {GELE}.
Unlike in the Lorenz equations, we observe limit tori, which are the multi-dimensional analogue of limit cycles, in the solutions of the DNS and {GELE} at high orders. 
Initial condition dependency in the {DNS} and Lorenz equations is also discussed.
\end{abstract}

\maketitle

\begin{quotation}
The Lorenz equations are a simplified nonlinear dynamical system derived from the two-dimensional Rayleigh-B\'enard convection problem. 
They have been one of the best-known examples in chaos theory due to the peculiar bifurcation and chaos behaviors. 
And they are often regarded as the minimal chaotic model for describing the convection system and, by extension, weather. 
Such an interpretation is sometimes challenged due to the simplifying restriction of considering only a few harmonics in the derivation.
This study loosens this restriction by considering additional high-order harmonics and derives a system we call {the generalized expansion of the Lorenz equations (GELE)}. 
{GELE} allows us to study how solutions transition from the classic Lorenz equations to high-order systems comparable to a two-dimensional Direct Numerical Simulation (DNS). 
This study also proposes mathematical formulations for a direct comparison between the Lorenz equations, {GELE}, and two-dimensional DNS as the system's order increases.
This work advances our understanding of the convection system by bridging the gap between the classic model of Lorenz and a more realistic convection system.
\end{quotation}

\section{\label{sec:Intro}Introduction}
The Rayleigh-B\'enard (RB) system is a canonical example of a flow convection system driven by the temperature difference $\Delta T$ between two boundaries in a plane horizontal fluid layer. 
When this condition of having higher temperature (i.e. $\Delta T>0$) and lower density at the bottom is maintained, such an unstable environment created by the thermal stratification can introduce a roll-type convection motion for a high enough $\Delta T$.
In more precise terms, the onset of convection motion happens when the nondimensional Rayleigh number $\mathrm{Ra}$, the ratio between buoyancy force and viscous force, is above its critical value $\mathrm{Ra}_{c}$. The critical Rayleigh number $\mathrm{Ra}_{c}$ depends on the boundary conditions and other system configurations.
As $\mathrm{Ra}$ increases further above $\mathrm{Ra}_{c}$ (i.e. $r=\mathrm{Ra}/\mathrm{Ra}_{c}\gg 1$), the RB system exhibits very rich dynamical behaviors such as instability, bifurcation, turbulence, chaos, intermittency, etc.
Due to its simple configuration despite the flow's complex behavior, the RB system has remained a popular research topic for over a century in diverse scientific disciplines including fluid mechanics, applied mathematics, and atmospheric science \citep[][]{Getling1998,Bodenschatz2000}.

In 1962, \citet[][]{Saltzman1962} further simplified the governing equations of the two-dimensional RB system into a highly truncated system of ordinary differential equations, which was cast as an initial value problem by applying the Fourier representations. 
The spectral analysis allows us to better understand the convection roll by considering it as the primary mode together with its nonlinear interactions with higher-order Fourier modes. 
Although \citet{Saltzman1962} was first to propose these nonlinear dynamical equations, its lowest order formulation by \citet{Lorenz1963} called the Lorenz equations is more widely recognised due to its association with Lorenz's discovery of deterministic chaos.

It is said that Lorenz had realized by chance that the finite predictability of weather might lie in nonlinearity of the governing systems in some fundamental sense.
In order to best illustrate the idea that even a simple deterministic system can exhibit sensitive initial-condition dependency and is therefore unpredictable, Lorenz settled on a system of three ordinary differential equations derived from the two-dimensional RB system, now known as the Lorenz equations. 
Being simple and deterministic, its derivation is still strongly rooted in the physics of thermal convection, following the Fourier-Galerkin method of approximating the governing equations for the two-dimensional RB system. 
As such, the Rayleigh number retains its relevance through the normalized Rayleigh number $r$, an important parameter controlling the onset of chaos in the Lorenz equations.
The butterfly-shaped Lorenz attractor \citep[][]{tucker1999lorenz} is arguably the most prominent image of chaos theory, the field which by mid 1980s morphed itself into some kind of a new scientific movement with profound and lasting influences across different disciplines \citep[][]{gleick1988chaos}.

More recently, efforts have been made to understand how nonlinear dynamical systems behave when the dimension of nonlinear dynamical systems increases. 
For instance, \citet{Shen2014} extended the Lorenz equations by incorporating two additional higher-order Fourier modes and studied their influence on the system. The nonlinear dynamical systems can also be extended by considering additional physical effects (e.g. rotation, scalar diffusion) in the governing equations\citep[][]{Stenflo1996,Park2016PS,Moon2019}. 
These extended systems exhibit somewhat different and sometimes new dynamical behaviors compared to the low-order Lorenz equations. For example, \citet{Felicio2018} demonstrated that a six-dimensional Lorenz-like system can even exhibit hyperchaos, (i.e. solutions with at least two positive Lyapunov exponents, which was not seen in the original Lorenz equations). 
For a systematic comparison between the classic Lorenz equations and the higher-order extensions, \citet{Moon2017} thoroughly investigated the dynamical behaviors and bifurcation structures of the extended systems obtained by considering higher-order harmonics at dimensions 5, 6, 8, 9, and 11 in wide ranges of parameters, which was later generalized \citep[][]{Moon2020} into explicit ODE expressions for $(3N)$- and $(3N+2)$-dimensional Lorenz systems for any positive integer $N$.

Two issues, however, remain unresolved in such analyses of the extensions at higher dimensions. 
First, as with all Lorenz and high-order Lorenz-like systems, it is not well-understood how much of the two-dimensional RB convection remains intact under the conversion into the Lorenz equations even at very high dimensions.
Conversely, it is also important to assess to what extent the many interesting nonlinear phenomena observed in the Lorenz equations are also found in the two-dimensional RB convection.
This study aims to address this issue by directly comparing the solutions of the Lorenz equations with results from a Direct Numerical Simulation (DNS) of the two-dimensional RB convection using the governing equations. 
There have been a number of DNS studies on the 2D RB convection \citep[][]{Stevens2011,Bao2017}, but most focus on instabilities and turbulence phenomena; explicit investigations about similarities and differences between the Lorenz equations and DNS have been rare still.
\citet{Paul2012} reported some bifurcation characteristics in the $r$ parameter space reminiscent of the Lorenz equations using the DNS. 
Nevertheless, a systematic and comparative investigation of the classic Lorenz equations and the DNS is still missing.

The second issue is pertinent to the way in which the dimension is raised in the previously investigated generalizations of the Lorenz equations \citep[]{Moon2017,Moon2020}, wherein the additionally incorporated higher-order harmonics are exclusively in the vertical direction of the thermal convection problem. 
These studies have not simultaneously considered horizontal higher-order harmonics and consequently the convection cells corresponding to very high harmonics in their generalizations may appear to have been vertically squeezed, which can lead to certain unnatural behaviors with regard to fluid convection. 
In this study, we newly formulate the {generalized expansion of the Lorenz equations (GELE)} by simultaneously considering higher-order harmonics in both the vertical and horizontal directions. 
{GELE} will serve as a link between the classic Lorenz equations and the DNS and will allow us a more complete investigation of the impact of higher-order harmonics on the various dynamical behaviors observed in the Lorenz equations.

The formulations of the equations for the DNS and {GELE} necessary for the systematic analysis are presented in Section \ref{sec:Prob}. 
Detailed descriptions on the governing equations, the modal amplitudes, energy relations, etc., are provided for the three different systems: the Lorenz equations, the DNS, and {GELE}. 
In Section \ref{sec:Numerical}, we demonstrate various numerical results; for instance, chaotic and equilibrium solutions, solution transition from the Lorenz equations to the DNS via variations of the order of {GELE}, periodic nature of the high-order systems, and initial-condition dependency. 
Finally in Section \ref{sec:Conclusion}, conclusions and discussion are given.

\section{\label{sec:Prob}Problem formulation}
\subsection{\label{sec:Prob_primitive}Primitive equations}
In the Cartesian coordinate $(x,z)$ where $x$ and $z$ are the streamwise (horizontal) and vertical coordinates, respectively, we consider the two-dimensional Navier-Stokes equations under the Boussinesq approximation together with the thermal diffusion equation as follows:
\begin{equation}
\label{eq:continuity}
\frac{\partial u}{\partial x}+\frac{\partial w}{\partial z}=0,
\end{equation}
\begin{equation}
\label{eq:x_mom}
\frac{\partial u}{\partial t}+u\frac{\partial u}{\partial x}+w\frac{\partial u}{\partial z}=-\frac{1}{\rho}_{0}\frac{\partial P}{\partial x}+\nu_{0}\nabla^{2}u,
\end{equation}
\begin{equation}
\label{eq:z_mom}
\frac{\partial w}{\partial t}+u\frac{\partial w}{\partial x}+w\frac{\partial w}{\partial z}=-\frac{1}{\rho_{0}}\frac{\partial P}{\partial z}-\frac{\Delta \bar{\rho}}{\rho_{0}} g+\nu_{0}\nabla^{2}w,
\end{equation}
\begin{equation}
\label{eq:energy}
\frac{\partial T}{\partial t}+u\frac{\partial T}{\partial x}+w\frac{\partial T}{\partial z}=\kappa_{0}\nabla^{2}T,
\end{equation}
where $u$ is the streamwise velocity, $w$ is the vertical velocity, $P$ is the pressure, $T$ is the temperature, $\Delta\bar{\rho}=\rho-\rho_{0}$ is the deviation of the density $\rho$ from the reference density $\rho_{0}$, $\nu_{0}$ is the reference kinematic viscosity, $\kappa_{0}$ is the thermal diffusivity, and $\nabla^{2}=\partial^{2}/\partial x^{2}+\partial^{2}/\partial z^{2}$ is the Laplacian operator. 
The reference values are computed from the properties at the bottom boundary $z=0$.
We assume that the density $\rho$ and the temperature $T$ satisfy a linear relation
\begin{equation}
\frac{\rho-\rho_{0}}{\rho_{0}}=-\epsilon_{0}\left(T-T_{0}\right),
\end{equation}
where $\epsilon_{0}$ is the thermal expansion coefficient and $T_{0}$ is the reference temperature. 
We assume that the temperature $T$ is given as
\begin{equation}
T=T_{0}-\frac{\Delta T}{H}z+\theta,
\end{equation}
where $\Delta T=T_{0}-T|_{z=H}>0$ is the temperature difference between $z=0$ and $z=H$ where $H$ is the domain height, and $\theta$ is the temperature perturbation. 
The pressure $P$ is assumed to be decomposed into $P=\mathcal{P}+p$ where $\mathcal{P}$ is the pressure satisfying the hydrostatic balance: $\partial\mathcal{P}/\partial z=-\rho_{0}\epsilon_{0}g\Delta T (z/H)$, and $p$ is the pressure perturbation.
Applying the above assumptions, we obtain the following set of equations:
\begin{equation}
\label{eq:continuity_2nd}
\frac{\partial u}{\partial x}+\frac{\partial w}{\partial z}=0,
\end{equation}
\begin{equation}
\label{eq:x_mom_2nd}
\frac{\partial u}{\partial t}+u\frac{\partial u}{\partial x}+w\frac{\partial u}{\partial z}=-\frac{1}{\rho}_{0}\frac{\partial p}{\partial x}+\nu_{0}\nabla^{2}u,
\end{equation}
\begin{equation}
\label{eq:z_mom_2nd}
\frac{\partial w}{\partial t}+u\frac{\partial w}{\partial x}+w\frac{\partial w}{\partial z}=-\frac{1}{\rho_{0}}\frac{\partial p}{\partial z}+\epsilon_{0}g\theta+\nu_{0}\nabla^{2}w,
\end{equation}
\begin{equation}
\label{eq:thermal_2nd}
\frac{\partial \theta}{\partial t}+u\frac{\partial \theta}{\partial x}+w\frac{\partial \theta}{\partial z}-\frac{\Delta T}{H}w=\kappa_{0}\nabla^{2}\theta.
\end{equation}
To analyze the system in a nondimensional form, we consider the reference time scale as $H^{2}/\kappa_{0}$, the length scale as $H$, the velocity scale as $\kappa_{0}/H$, the pressure scale as $\rho_{0}\kappa_{0}^{2}/H^{2}$, and the temperature scale $\Delta T$. 
Then the nondimensional equations read
\begin{equation}
\label{eq:continuity_dimless}
\frac{\partial u}{\partial x}+\frac{\partial w}{\partial z}=0,
\end{equation}
\begin{equation}
\label{eq:x_mom_dimless}
\frac{\partial u}{\partial t}+u\frac{\partial u}{\partial x}+w\frac{\partial u}{\partial z}=-\frac{\partial p}{\partial x}+\sigma\nabla^{2}u,
\end{equation}
\begin{equation}
\label{eq:z_mom_dimeless}
\frac{\partial w}{\partial t}+u\frac{\partial w}{\partial x}+w\frac{\partial w}{\partial z}=-\frac{\partial p}{\partial z}+\sigma\mathrm{Ra}\theta+\sigma\nabla^{2}w,
\end{equation}
\begin{equation}
\label{eq:thermal_dimless}
\frac{\partial \theta}{\partial t}+u\frac{\partial \theta}{\partial x}+w\frac{\partial \theta}{\partial z}-w=\nabla^{2}\theta,
\end{equation}
where $\sigma=\nu_{0}/\kappa_{0}$ is the Prandtl number and $\mathrm{Ra}=\epsilon_{0}gH^{3}\Delta T /\kappa_{0}\nu_{0}$ is the Rayleigh number. 
Note that the variables $(u,w,p,\theta)$ are now dimensionless.
The set of equations (\ref{eq:continuity_dimless})--(\ref{eq:thermal_dimless}) can be further simplified if we consider the streamfunction $\psi$ that satisfies
\begin{equation}
u=-\frac{\partial\psi}{\partial z},~w=\frac{\partial\psi}{\partial x}.
\end{equation}
The simplified set of equations for $\psi$ and $\theta$ becomes
\begin{equation}
\label{eq:psi}
\frac{\partial }{\partial t}\nabla^{2}\psi=\frac{\partial\psi}{\partial z}\frac{\partial\nabla^{2}\psi}{\partial x}-\frac{\partial\psi}{\partial x}\frac{\partial \nabla^{2}\psi}{\partial z}+\sigma\nabla^{4}\psi+\sigma \mathrm{Ra}\frac{\partial \theta}{\partial x},
\end{equation}
\begin{equation}
\label{eq:theta}
\frac{\partial\theta}{\partial t}=\frac{\partial\psi}{\partial z}\frac{\partial\theta}{\partial x}-\frac{\partial\psi}{\partial x}\frac{\partial \theta}{\partial z}+\nabla^{2}\theta+\frac{\partial\psi}{\partial x},
\end{equation}
(see also, \citet{Saltzman1962}).

We solve the equations (\ref{eq:psi})--(\ref{eq:theta}) by imposing the boundary conditions such that variables $\psi$ and $\theta$ are periodic in the $x$-direction:
\begin{equation}
\label{eq:bc_x}
\psi(x=0,z)=\psi(x=l_{x},z),~
\theta(x=0,z)=\theta(x=l_{x},z),
\end{equation}
where $l_{x}$ is the streamwise domain length, while we consider in the $z$-direction the following boundary conditions
\begin{equation}
\label{eq:bc_z}
\psi=\theta=\frac{\partial^{2}\psi}{\partial z^{2}}=0,
\end{equation}
at $z=0$ and $z=1$.
The equations (\ref{eq:psi})--(\ref{eq:theta}) in the physical space $(x,z)$ as well as the boundary conditions (\ref{eq:bc_x})--(\ref{eq:bc_z}) will be used in the two-dimensional DNS. 
And we will describe in the last subsection the numerical methods for performing the two-dimensional DNS. 

\subsection{Relation between DNS and Lorenz formulations}
For the derivation of the classic Lorenz equations, we consider the following transformations 
\begin{eqnarray}
\label{eq:Lorenz_transform}
\psi(x,z,t)&=&X(t)\frac{\sqrt{2}(\alpha^{2}+\beta^{2})}{\alpha\beta}\sin(\alpha x)\sin(\beta z),\nonumber\\
\theta(x,z,t)&=&Y(t)\frac{\sqrt{2}(\alpha^{2}+\beta^{2})^{3}}{\alpha^{2}\beta\mathrm{Ra}}\cos(\alpha x)\sin(\beta z)\nonumber\\
&-&Z(t)\frac{(\alpha^{2}+\beta^{2})^{3}}{\alpha^{2}\beta\mathrm{Ra}}\sin(2\beta z), 
\end{eqnarray}
where $(X,Y,Z)$ are the time-dependent amplitudes, $\alpha=2\pi/l_{x}$ is the streamwise wavenumber, and $\beta=\pi$ is the vertical wavenumber. 
Note that the above transformations truncate off other high-order harmonics in the $x$- and $z$-directions. 
Using (\ref{eq:Lorenz_transform}) and neglecting high-order nonlinear interactions as such, we derive the Lorenz equations:
\begin{eqnarray}
\label{eq:Lorenz}
\frac{\mathrm{d}X}{\mathrm{d}\tau}&=&\sigma(Y-X),\nonumber\\
\frac{\mathrm{d}Y}{\mathrm{d}\tau}&=&rX-Y-XZ,\nonumber\\
\frac{\mathrm{d}Z}{\mathrm{d}\tau}&=&XY-bZ,
\end{eqnarray}
where $\tau=(\alpha^{2}+\beta^{2})t$ is the rescaled time, $r=\mathrm{Ra}/\mathrm{Ra}_{c}$ is the normalized Rayleigh number (i.e. the ratio between the Rayleigh number and the critical Rayleigh number $\mathrm{Ra}_{c}=(\alpha^{2}+\beta^{2})^{3}/\alpha^{2}$), and $b=4\beta^{2}/(\alpha^{2}+\beta^{2})$ is the geometrical parameter.

Once we solve the Lorenz equations (\ref{eq:Lorenz}), we can recover the Lorenz-based physical solutions {$\psi^{(\mathrm{Lo})}(x,z)$ and $\theta^{(\mathrm{Lo})}(x,z)$} by using the backward transformations (\ref{eq:Lorenz_transform}).
Since nonlinear interactions among high-order harmonics are ignored, {$\psi^{(\mathrm{Lo})}$ and $\theta^{(\mathrm{Lo})}$} are different from those $\psi$ and $\theta$ obtained from the DNS.
To quantify the differences more systematically, we compute the DNS-based amplitudes $(X^{(\mathrm{D})}, Y^{(\mathrm{D})}, Z^{(\mathrm{D})})$ as follows:
\begin{eqnarray}
\label{eq:DNS_XYZ}
X^{(\mathrm{D})}&=&\frac{\sqrt{2}\alpha^{2}\beta}{\pi(\alpha^{2}+\beta^{2})}\int_{0}^{l_{x}}\int_{0}^{1}\psi\sin(\alpha x)\sin(\beta z)\mathrm{d}z\mathrm{d}x,\nonumber\\
Y^{(\mathrm{D})}&=&\frac{\sqrt{2}\alpha^{3}\beta\mathrm{Ra}}{\pi(\alpha^{2}+\beta^{2})^{3}}\int_{0}^{l_{x}}\int_{0}^{1}\theta\cos(\alpha x)\sin(\beta z)\mathrm{d}z\mathrm{d}x,\nonumber\\
Z^{(\mathrm{D})}&=&\frac{-\alpha^{3}\beta\mathrm{Ra}}{\pi(\alpha^{2}+\beta^{2})^{3}}\int_{0}^{l_{x}}\int_{0}^{1}\theta\sin(2\beta z)\mathrm{d}z\mathrm{d}x,
\end{eqnarray}
where $\psi(x,z)$ and $\theta(x,z)$ in (\ref{eq:DNS_XYZ}) are the variables computed from the DNS. 
Note that the DNS-based amplitudes $(X^{(\mathrm{D})},Y^{(\mathrm{D})},Z^{(\mathrm{D})})$ are obtained by integrations over the domain length in the vertical direction $z$ and one wavelength in the streamwise direction $x$.

\subsection{\label{sec:Prob_DNS}Spectral formulation for generalized nonlinear dynamical system}
In this study, we assume that the solution is spatially periodic in the $x$-direction and bounded in the $z$-direction as a way to allow the Fourier representations \citep[][]{Saltzman1962}. 
This consideration allows us to express the physical solution $\psi$ and $\theta$ in the spectral form. 
First, we consider the spatial periodicity in the $x$-direction by expressing $\psi$ and $\theta$ as
\begin{equation}
\label{eq:mode_dns}
\left(
\begin{array}{c}
\psi(x,z,t)\\
\theta(x,z,t)
\end{array}
\right)=\sum_{l=-L}^{L}
\left(
\begin{array}{c}
\tilde{\psi}_{l}(z,t)\\
\tilde{\theta}_{l}(z,t)
\end{array}
\right)\exp(\mathrm{i}\alpha_{l}x),
\end{equation}
where $l$ is the mode number, $L$ is the largest mode number we consider for the streamwise spectral modes, $\tilde{\psi}_{l}(z,t)$ and $\tilde{\theta}_{l}(z,t)$ are the mode shapes of $\psi$ and $\theta$, respectively, $\mathrm{i}=\sqrt{-1}$, and $\alpha_{l}=l\alpha$ is the streamwise wavenumber of the mode $l$.
Since $\psi$ and $\theta$ are real, the complex-conjugate modal relations $\tilde{\psi}_{-l}=\tilde{\psi}_{l}^{*}$ and $\tilde{\theta}_{-l}=\tilde{\theta}_{l}^{*}$ (where $*$ denotes the complex conjugate) must be satisfied for $l\geq1$, while $\tilde{\psi}_{0}$ and $\tilde{\theta}_{0}$ must be real.
For each mode $l$, we express the equations (\ref{eq:psi}) and (\ref{eq:theta}) in the modal form as 
\begin{equation}
\label{eq:psi_spec_x}
\frac{\partial }{\partial t}\tilde{\nabla}^{2}_{l}\tilde{\psi}_{l}=\sigma\tilde{\nabla}^{4}_{l}\tilde{\psi}_{l}+\mathrm{i}\alpha_{l}\sigma \mathrm{Ra}\tilde{\theta}_{l}+\tilde{N}^{\psi}_{l},
\end{equation}
\begin{equation}
\label{eq:theta_spec_x}
\frac{\partial\tilde{\theta}_{l}}{\partial t}=\tilde{\nabla}^{2}_{l}\tilde{\theta}_{l}+\mathrm{i}\alpha_{l}\tilde{\psi}_{l}+\tilde{N}^{\theta}_{l},
\end{equation}
where $\tilde{\nabla}^{2}_{l}=\frac{\partial^{2}}{\partial z^{2}}-\alpha_{l}^{2}$, and $\tilde{N}^{\psi}_{l}$ and $\tilde{N}^{\theta}_{l}$ are the convolution terms:
\begin{eqnarray}
\label{eq:nonlinear_dns}
\tilde{N}^{\psi}_{l}&=&\sum_{j=-L}^{L}\mathrm{i}\alpha_{j}\left[\tilde{\nabla}^{2}_{j}\tilde{\psi}_{j}\frac{\partial \tilde{\psi}_{l-j}}{\partial z}-\tilde{\psi}_{j}\frac{\partial}{\partial z}\left(\tilde{\nabla}^{2}_{l-j}\tilde{\psi}_{l-j}\right)\right],\nonumber\\
\tilde{N}^{\theta}_{l}&=&\sum_{j=-L}^{L}\mathrm{i}\alpha_{j}\left(\tilde{\theta}_{j}\frac{\partial\tilde{\psi}_{l-j}}{\partial z}-\tilde{\psi}_{j}\frac{\partial\tilde{\theta}_{l-j}}{\partial z}\right),
\end{eqnarray}
which are related to the nonlinear terms in (\ref{eq:psi}) and (\ref{eq:theta}).
Note that, in the spectral transformation (\ref{eq:mode_dns}) and the nonlinear convolution (\ref{eq:nonlinear_dns}), high-order harmonics ($|l|>L$) generated by nonlinear interactions of low-order harmonics ($|l|\leq L$) are ignored. 
In principle, the spectral solution in the limit $L\rightarrow\infty$ will recover the DNS solution in the physical space $(x,z)$. 
On the other hand, if $L=1$, the spectral solution can match the Lorenz solution when low-order harmonics in the $z$-direction are considered. 
The mode number limit $L$ is, therefore, an important control parameter that allows us to study the transition from the Lorenz equations to the DNS. 

The ansatz (\ref{eq:mode_dns}) is spectral only in the $x$-direction but we can further expand the mode shapes $\tilde{\psi}$ and $\tilde{\theta}$ using the sinuous series in the $z$-direction as follows:
\begin{equation}
\label{eq:modal_xz}
\left(
\begin{array}{c}
\tilde{\psi}_{l}(z,t)\\
\tilde{\theta}_{l}(z,t)
\end{array}
\right)=\sum_{m=0}^{M}
\left(
\begin{array}{c}
\hat{\psi}_{lm}(t)\\
\hat{\theta}_{lm}(t)
\end{array}
\right)\sin(\beta_{m} z),
\end{equation}
where $\hat{\psi}_{lm}$ and $\hat{\theta}_{lm}$ are the time-dependent mode amplitudes, $m$ is the mode number in the $z$-direction, $M$ is the largest mode number we consider for the vertical spectral modes, and $\beta_{m}=m\beta$ is the vertical wavenumber of the mode $m$.
Note that the sinuous series with $\sin(\beta_{m}z)$ satisfies the boundary conditions at $z=0$ and 1 for any $m$.
Applying the expansion (\ref{eq:modal_xz}) to the equations (\ref{eq:psi_spec_x})-(\ref{eq:theta_spec_x}) leads to the following equations of the {generalized expansion of the Lorenz equations}:
\begin{equation}
\label{eq:modal_psi}
-\left(\alpha_{l}^{2}+\beta_{m}^{2}\right)\frac{\mathrm{d} \hat{\psi}_{lm}}{\mathrm{d} t}=\sigma \left(\alpha_{l}^{2}+\beta_{m}^{2}\right)^{2}\hat{\psi}_{lm}+\mathrm{i}\alpha_{l}\sigma \mathrm{Ra}\hat{\theta}_{lm}+\hat{N}^{\psi}_{lm},
\end{equation}
\begin{equation}
\label{eq:modal_theta}
\frac{\mathrm{d} \hat{\theta}_{lm}}{\mathrm{d} t}=-\left(\alpha_{l}^{2}+\beta_{m}^{2}\right)\hat{\theta}_{lm}+\mathrm{i}\alpha_{l}\hat{\psi}_{lm}+\hat{N}_{lm}^{\theta},
\end{equation}
where $\hat{N}^{\psi}_{lm}$ and $\hat{N}^{\theta}_{lm}$ are the convolution terms derived from the nonlinear terms $\tilde{N}_{l}^{\psi}$ and $\tilde{N}_{l}^{\theta}$ (see Appendix \ref{sec:Appendix_nonlinear} for more details).

The practicality of the {GELE} above is in that the equations (\ref{eq:modal_psi})--(\ref{eq:modal_theta}) can produce either the DNS solutions or the Lorenz solutions depending on the choice of $L$ and $M$.
For instance, {GELE} can be simplified into the Lorenz equations when we consider $L=1$ and $M=2$ and when proper initial conditions are imposed such that initial mode amplitudes except $\Im(\hat{\psi}_{11})$, $\Re(\hat{\theta}_{11})$ and $\hat{\theta}_{02}$ are zero (i.e. $\Re(\hat{\psi}_{11})=\Im(\hat{\theta}_{11})=0$, $\hat{\psi}_{01}=\hat{\psi}_{02}=\hat{\psi}_{12}=\hat{\theta}_{01}=\hat{\theta}_{12}=0$ where $\Re$ and $\Im$ denote the real and imaginary parts, respectively).
As similarly derived for the DNS-based amplitudes in (\ref{eq:DNS_XYZ}), the {GELE}-based amplitudes $X^{(\mathrm{G})}$, $Y^{(\mathrm{G})}$ and $Z^{(\mathrm{G})}$ can be computed from the following relations:
\begin{eqnarray}
X^{(\mathrm{G})}(t)&=&-\frac{\sqrt{2}\alpha\beta}{\left(\alpha^{2}+\beta^{2}\right)}\Im\left[\hat{\psi}_{11}(t)\right],\nonumber\\
Y^{(\mathrm{G})}(t)&=&\frac{\sqrt{2}\alpha^{2}\beta\mathrm{Ra}}{\left(\alpha^{2}+\beta^{2}\right)^{3}}\Re\left[\hat{\theta}_{11}(t)\right],\nonumber\\
Z^{(\mathrm{G})}(t)&=&-\frac{\alpha^{2}\beta\mathrm{Ra}}{\left(\alpha^{2}+\beta^{2}\right)^{3}}\hat{\theta}_{02}(t).
\end{eqnarray}
If we consider $M>2$ and $L=1$, we recover the high-order Lorenz equations \citep[][]{Moon2017,Moon2020}. 
And we can also reproduce the results of the DNS mathematically in the limits $L\rightarrow\infty$ and $M\rightarrow\infty$ (in practice, when $L$ and $M$ are sufficiently large).
Furthermore, the mode amplitudes in {GELE} can be directly compared with those from the DNS if we consider the DNS-based mode amplitudes $\hat{\psi}_{lm}^{(\mathrm{D})}$ and $\hat{\theta}_{lm}^{(\mathrm{D})}$ obtained from the following relations:
\begin{eqnarray}
\hat{\psi}_{lm}^{(\mathrm{D})}&=&\frac{\alpha}{\pi}\int_{0}^{l_{x}}\int_{0}^{1}\psi\sin(\beta_{m}z)\exp(-\mathrm{i}\alpha_{l}x)\mathrm{d}z\mathrm{d}x,\nonumber\\
\hat{\theta}_{lm}^{(\mathrm{D})}&=&\frac{\alpha}{\pi}\int_{0}^{l_{x}}\int_{0}^{1}\theta\sin(\beta_{m}z)\exp(-\mathrm{i}\alpha_{l}x)\mathrm{d}z\mathrm{d}x.
\end{eqnarray}

\subsection{Dissipative system and energy relations}
By taking the divergence, we can check whether {GELE} is dissipative \citep[][]{Lorenz1963}. 
Applying the partial derivatives of the equations (\ref{eq:modal_psi}) and (\ref{eq:modal_theta}) with respect to $\hat{\psi}_{lm}$ and $\hat{\theta}_{lm}$, we have
\begin{eqnarray}
\label{eq:volume_contraction}
&&\sum_{l=-L}^{L}\sum_{m=0}^{M}\left[\frac{\partial}{\partial\hat{\psi}_{lm}}\left(\frac{\mathrm{d}\hat{\psi}_{lm}}{\mathrm{d}t}\right)+\frac{\partial}{\partial\hat{\theta}_{lm}}\left(\frac{\mathrm{d}\hat{\theta}_{lm}}{\mathrm{d}t}\right)\right]\nonumber\\
&&=-(\sigma+1)\sum_{l=-L}^{L}\sum_{m=0}^{M}\left(\alpha_{l}^{2}+\beta_{m}^{2}\right).
\end{eqnarray}
We clearly see that the right-hand-side term is always negative, which implies that the system is dissipative. 
As similarly pointed out by \citet{Moon2017}, the right-hand-side term of (\ref{eq:volume_contraction}) becomes largely negative and the volume contraction occurs at a faster rate when the limits of the system's order $L$ and $M$ increase.

It is also important to define the total energy $\mathrm{E}_{\mathrm{T}}$ which is the sum of the kinetic energy $\mathrm{E}_{\mathrm{K}}$ and potential energy $\mathrm{E}_{\mathrm{P}}$ (i.e. $\mathrm{E}_{\mathrm{T}}=\mathrm{E}_{\mathrm{K}}+\mathrm{E}_{\mathrm{P}}$), where these energies can be defined in dimensionless forms,
\begin{equation}
\label{eq:def_energies}
\mathrm{E}_{\mathrm{K}}=\int_{0}^{1}\int_{0}^{l_{x}}\frac{1}{2}\left(u^{2}+w^{2}\right)\mathrm{d}x\mathrm{d}z,~~
\mathrm{E}_{\mathrm{P}}=\int_{0}^{1}\int_{0}^{l_{x}}(-\sigma\mathrm{Ra}z)\theta \mathrm{d}x\mathrm{d}z.
\end{equation}
We note that the definition of $\mathrm{E}_{\mathrm{P}}$ above is different from that of \citet{Saltzman1962}, which is based on the square of the temperature perturbation. 
After manipulating the equations (\ref{eq:continuity_dimless})--(\ref{eq:thermal_dimless}) and considering the boundary conditions, the temporal evolution of the total energy can be written as follows:
\begin{equation}
\label{eq:def_energy_evolution}
\frac{\partial\mathrm{E}_{\mathrm{T}}}{\partial t}=\int_{0}^{1}\int_{0}^{l_{x}}\left(u\frac{\partial u}{\partial t}+w\frac{\partial w}{\partial t}-\sigma\mathrm{Ra}z\frac{\partial\theta}{\partial t}\right)\mathrm{d}x\mathrm{d}z=\mathcal{Q}+\mathcal{V},
\end{equation}
where $\mathcal{Q}$ is the temporal energy rate due to the thermal conduction occurring at the boundary $z=1$:
\begin{equation}
\label{eq:def_thermal_conduction}
\mathcal{Q}=-\sigma\mathrm{Ra}\int_{0}^{l_{x}}\left.z\frac{\partial\theta}{\partial z}\right|_{z=1}\mathrm{d}x,
\end{equation}
and $\mathcal{V}$ is the temporal energy rate due to the viscous dissipation:
\begin{eqnarray}
\label{eq:def_viscous_dissipation}
\mathcal{V}&=&-\sigma\int_{0}^{1}\int_{0}^{l_{x}}\left[\left(\frac{\partial u}{\partial x}\right)^{2}\right.+\left(\frac{\partial u}{\partial z}\right)^{2}\nonumber\\
&&+\left.\left(\frac{\partial w}{\partial x}\right)^{2}+\left(\frac{\partial w}{\partial z}\right)^{2}\right]\mathrm{d}x\mathrm{d}z.
\end{eqnarray}
It is important to note that $\mathcal{V}$ is always negative thus the viscous dissipation is responsible for the loss of the total energy, while $\mathcal{Q}$ can be positive or negative depending on the sign of the temperature gradient $\partial\theta/\partial z$ at $z=1$. 

If we use the spectral formulation (\ref{eq:modal_xz}), we can further simplify the energy expressions without integrations; for instance, we have the kinetic and potential energies
\begin{eqnarray}
\label{eq:energy_spectral}
\mathrm{E}_{\mathrm{K}}&=&\sum_{l=-L}^{L}\sum_{m=0}^{M}\frac{\pi\left(\alpha_{l}^{2}+\beta_{m}^{2}\right)}{2\alpha}|\hat{\psi}_{lm}|^{2},\nonumber\\
\mathrm{E}_{\mathrm{P}}&=&\sigma\mathrm{Ra}\sum_{m=1}^{M}\frac{2\pi\cos(\beta_{m})}{\alpha\beta_{m}}\hat{\theta}_{0m}.
\end{eqnarray}
Note that only the temperature modes $\hat{\theta}_{lm}$ with $l=0$ contribute to the potential energy since the integration in the $x$-direction in (\ref{eq:def_energies}) suppresses the contribution from the periodic modes $\hat{\theta}_{lm}$ of $l>0$.
The energy rates can be re-expressed as follows:
\begin{eqnarray}
\label{eq:energy_rate_spectral}
\mathcal{V}&=&-\sigma\sum_{l=-L}^{L}\sum_{m=0}^{M}\frac{\pi\left(\alpha_{l}^{2}+\beta_{m}^{2}\right)^{2}}{\alpha}|\hat{\psi}_{lm}|^{2},\nonumber\\
\mathcal{D}&=&-\sigma\mathrm{Ra}\sum_{m=1}^{M}\frac{2\pi\beta_{m}\cos(\beta_{m})}{\alpha}\hat{\theta}_{0m}.
\end{eqnarray}

\subsection{Numerical methods}
Considering the boundary conditions (\ref{eq:bc_x}) and (\ref{eq:bc_z}), we use the Chebyshev spectral method in the $z$-direction and the Fourier spectral method in the $x$-direction for numerical discretizations in the two-dimensional DNS \citep[][]{Weideman2000,Antkowiak2005,Park2017}. 
For the time stepping, we consider the implicit Euler method on the linear terms and the Adams-Bashforth scheme for the nonlinear terms \citep[][]{KMM1987}.
Direct numerical simulations in the physical space $(x,z)$ use an appropriate number of collocation points between 80 and 200 in both $x$- and $z$-directions and the time step $\Delta t$ between $10^{-6}$ and $10^{-4}$ in order to meet the Courant-Friedrichs-Lewy (CFL) condition for numerical stability in our parameter ranges of interest. 
When time-stepping {GELE} and the Lorenz equations, we also consider the implicit Euler method on the linear operator while the nonlinear terms are solved explicitly with the forward Euler method. 
For all results presented in this paper, some parameters such as $\sigma=10$ and $b=8/3$ are fixed (i.e. $\alpha=\pi/\sqrt{2}$ and $\beta=\pi$, the parameters that give $\mathrm{Ra}_{c}=27\pi^{4}/4$).
We only vary the parameters $r$, $L$ and $M$ as control parameters to elucidate the similarities and differences between the DNS, {GELE}, and the Lorenz equations.

In principle, a variety of types of initial conditions are available for numerical computation. 
For instance, we can impose Lorenz-like initial conditions where all the variables except $(X,Y,Z)$ are zero.
The Lorenz-like initial conditions in modal amplitudes can be converted into the DNS initial conditions as $\psi(x,z,0)=2|\hat{\psi}_{11}(0)|\sin(\alpha x)\sin(\beta z)$ and $\theta(x,z,0)=2|\hat{\theta}_{11}(0)|\cos(\alpha x)\sin(\beta z)+\hat{\theta}_{02}\sin(2\beta z)$. 
Although we can also impose various other kinds of initial conditions (e.g. non-zero higher harmonics where $\hat{\psi}_{lm}(0)\neq0$ or $\hat{\theta}_{lm}(0)\neq0$ or random initial conditions with random profiles of $\psi(x,z,0)$ and $\theta(x,z,0)$), we will mostly focus on the cases computed using the Lorenz-like conditions, and the initial condition sensitivity with random initial conditions will be discussed briefly.

\section{Numerical results}
\label{sec:Numerical}
We consider the regime $r>1$ (i.e. $\mathrm{Ra}>\mathrm{Ra}_{c}$), where the two-dimensional convection system is linearly unstable. 
As $r$ is increased from 1, we will investigate how dynamical behaviors such as bifurcation, nonlinear equilibration, chaos, or periodic attractors, all of which are only observable in the unstable regime and vary with the system orders $L$ and $M$. 
Note that when we say a regime is \emph{stable}, we refer to stability of the convection system not the stability of attractors.

\subsection{Chaotic and equilibrium states in the unstable regime}
\begin{figure*}
\includegraphics[height=6cm]{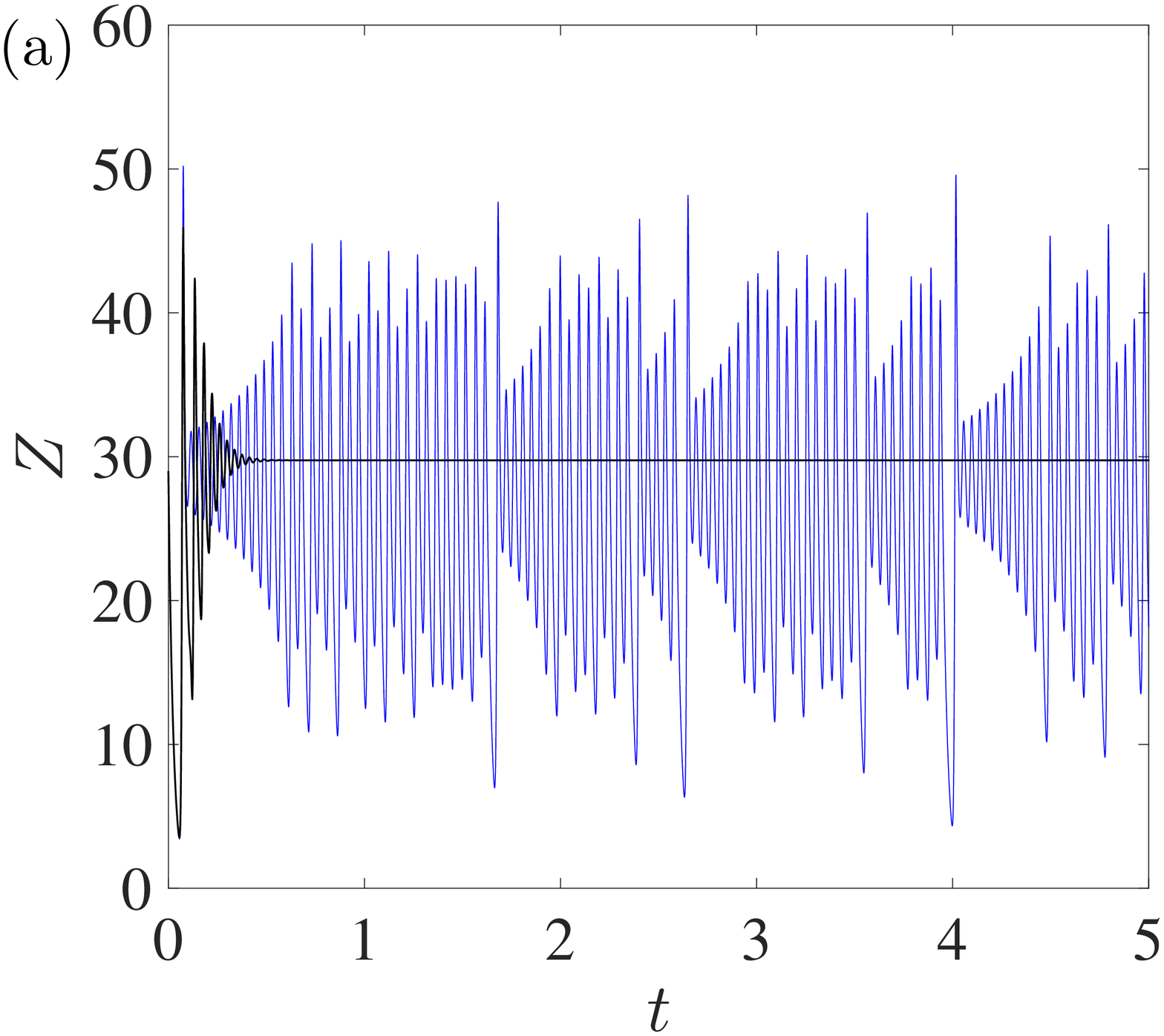}
\includegraphics[height=6cm]{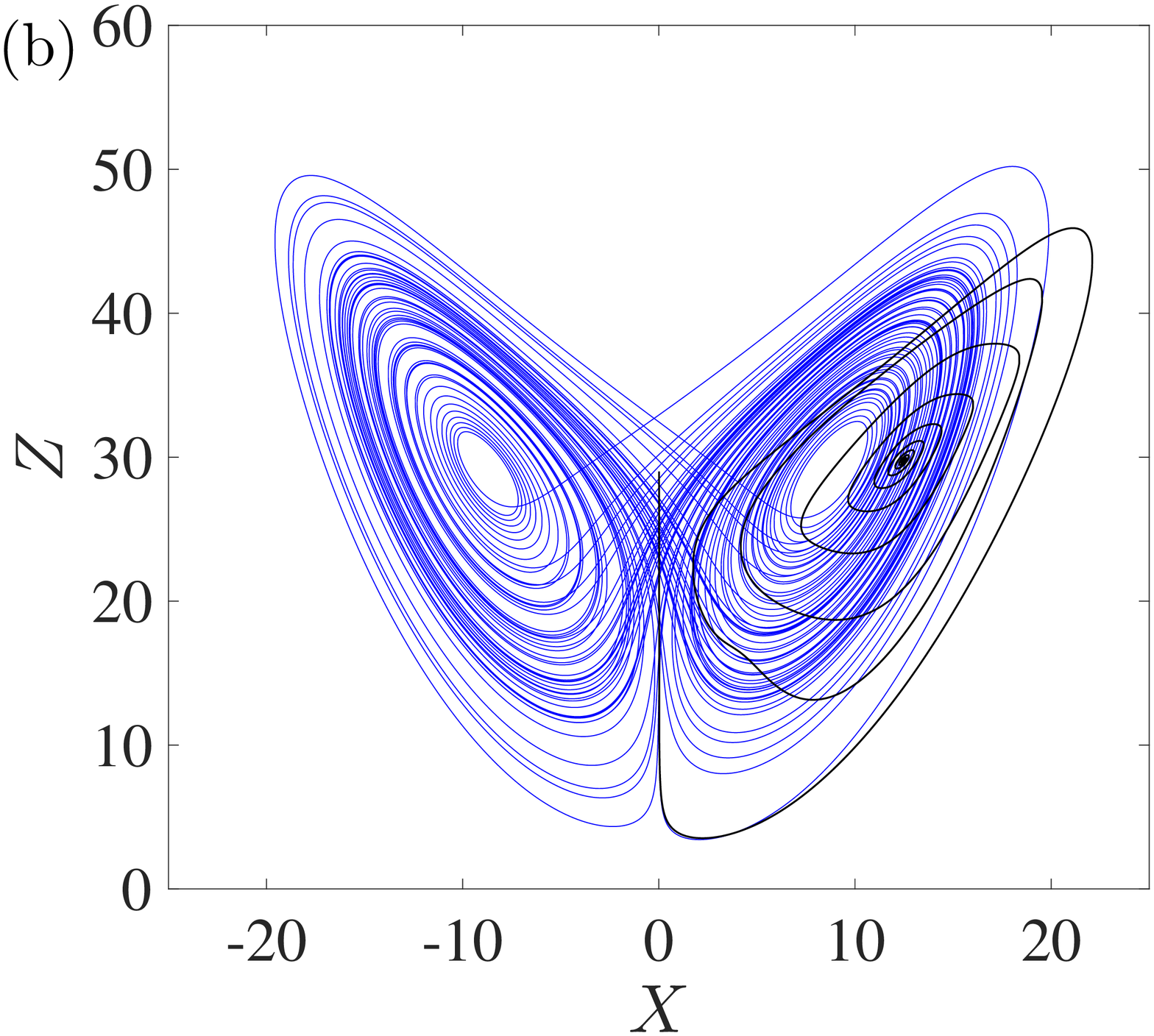}
\includegraphics[height=6cm]{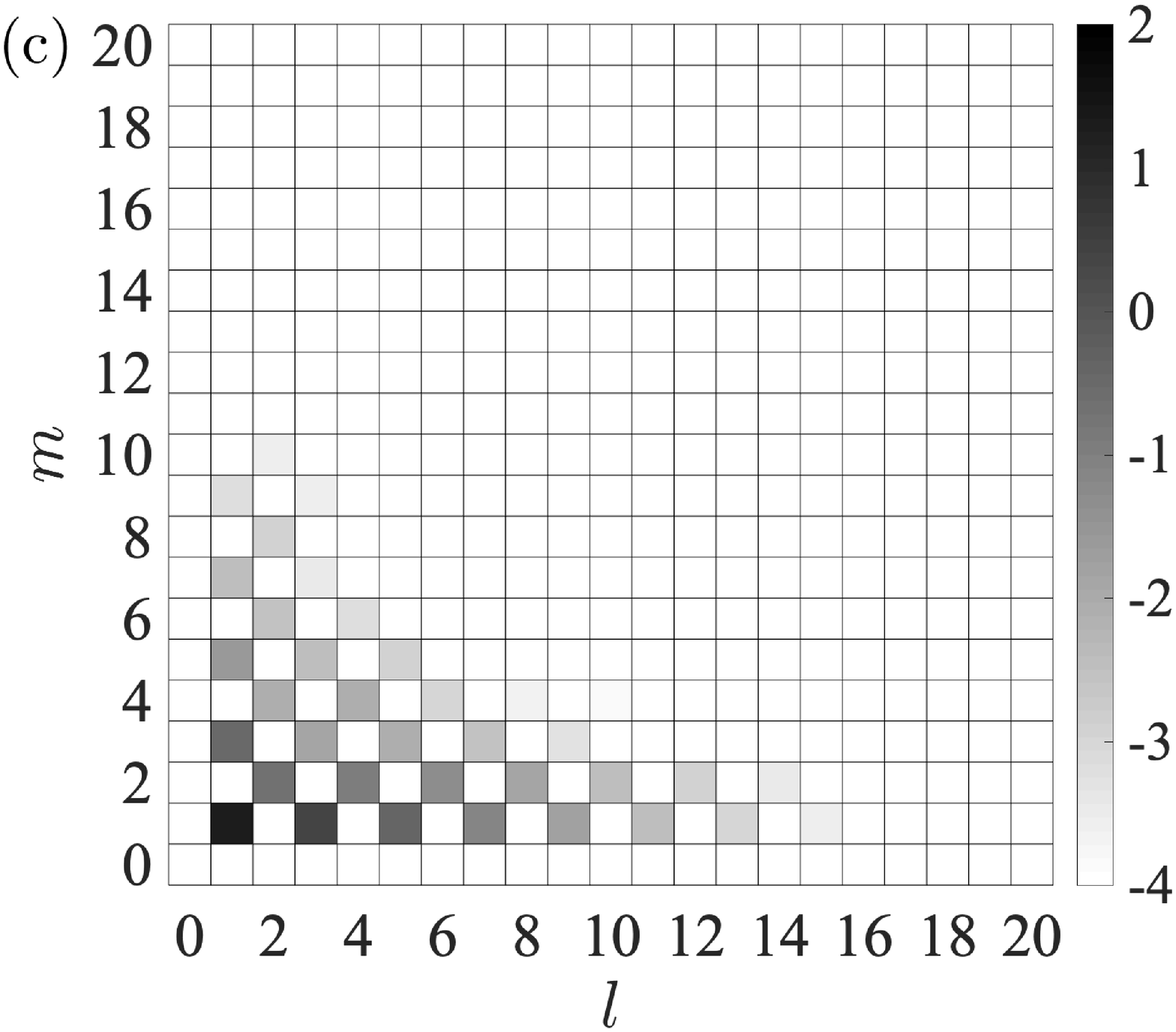}
\includegraphics[height=6cm]{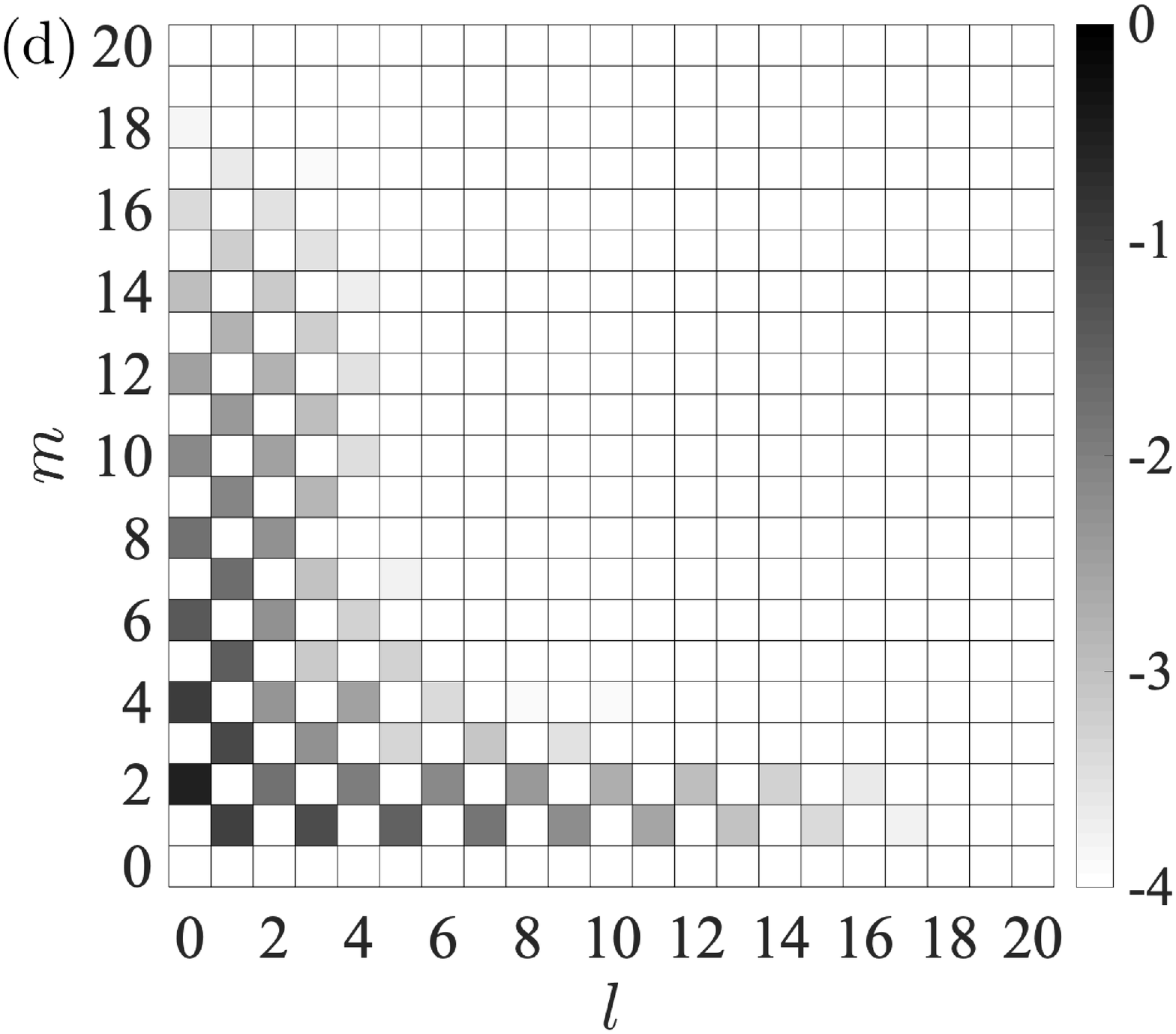}
\caption{\label{fig:case_r30} 
(a) Variable $Z$ versus time $t$ for the Lorenz solution (blue solid line) and the DNS solution (black solid line) at $r=30$. 
(b) Trajectories on the $(X,Z)$-plane of the Lorenz (blue) and DNS (black) solutions. 
(c,d) Amplitude distributions of the DNS solution: (c) $\log_{10}|\hat{\psi}_{lm}|$ and (d) $\log_{10}|\hat{\theta}_{lm}|$ in the parameter space of mode numbers $(l,m)$ at $t=5$. 
}
\end{figure*}
In this subsection, we fix $r=30$, a representative value at which we can observe the chaotic attractor in the classic Lorenz equations.
In Fig.~\ref{fig:case_r30}(a), we plot the amplitude $Z$ versus time $t$ and compare $Z(t)$ of the Lorenz equations with $Z^{(\mathrm{D})}(t)$ obtained from the DNS when the Lorenz-like initial condition $(X,Y,Z)=(0.01,0,r-1)$ is imposed on both the DNS and Lorenz equations. 
In fact, the temperature perturbation with $Y=0$ and $Z>0$ yields a stable solution when $X=0$, since the corresponding temperature solution in the physical space: $\theta(x,z,0)=\hat{\theta}_{02}(0)\sin(2\beta z)$ with $\hat{\theta}_{02}(0)<0$ implies that the temperature perturbation is stably stratified (i.e. $\theta$ is positive and the fluid density is lighter in the upper region $0.5<z<1$ while $\theta$ is negative and the fluid density is heavier in the lower region $0<z<0.5$).
However, we impose $X=0.01$ at $t=0$ to have a small-amplitude streamfunction perturbation, which has a roll shape and can cause instability.
Figure~\ref{fig:case_r30}(a) shows that there is a short transient period from $t=0$ where variable $Z$ decreases when $X$ is very small.
In this transient period, the DNS amplitude $Z^{(\mathrm{D})}(t)$ matches the Lorenz amplitude $Z(t)$, but afterwards $Z$ increases as $X$ is amplified and we see an oscillatory behavior of $Z$ in time $t$. 
A clear difference between the Lorenz equations and the DNS is now such that the Lorenz amplitude $Z$ becomes chaotic after the transient oscillatory period, while the DNS amplitude $Z^{(\mathrm{D})}$ reaches an equilibrium and converges to $Z^{(\mathrm{D})}\simeq29.75$ as $t$ increases.
These different dynamical behaviors can also be clearly distinguished in Fig.~\ref{fig:case_r30}(b), where the Lorenz solution exhibits a chaotic attractor on the $(X,Z)$-plane while the DNS solution moves along a spiral that converges to a fixed solution $(X^{(\mathrm{D})},Z^{(\mathrm{D})})\simeq(12.46,29.75)$.
We note that this DNS fixed solution is close to but is still different from the fixed point solution of the Lorenz equations: $(X,Z)|_{\mathrm{fixed}}=(\sqrt{b(r-1)},r-1)\simeq(8.79,29)$.
For variable $Y$, the DNS solution converges to $Y^{(\mathrm{D})}\simeq12.46$, a value still different from that of the fixed point solution $Y_{\mathrm{fixed}}=\sqrt{b(r-1)}\simeq8.79$ for the Lorenz equations.

The difference between the Lorenz and DNS solutions results from the fact that the DNS allows nonlinear interactions among higher-order modes. 
To see more clearly how the high-order nonlinear interactions occur in the DNS, we plot in Fig.~\ref{fig:case_r30}(c,d) the log-scale absolute values of the amplitudes $\hat{\psi}_{lm}$ and $\hat{\theta}_{lm}$ in the mode number space $(l,m)$ at $t=5$. 
Note that we only need to display the mode number space for non-negative $l\geq0$ due to the symmetries $\hat{\psi}_{(-l)m}^{*}=\hat{\psi}_{lm}$ and $\hat{\theta}_{(-l)m}^{*}=\hat{\theta}_{lm}$. 
The initial amplitudes we impose at $t=0$ are $X=0.01$ and $Z=r-1=29$ (i.e. $\hat{\psi}_{11}=-0.015\mathrm{i}$ and $\hat{\theta}_{02}=-0.3077$), while other variables are zero.
On the one hand, the Lorenz equations only allow nonlinear interactions between $\hat{\psi}_{11}$, $\hat{\theta}_{02}$, and $\hat{\theta}_{11}$. 
If we plot the amplitudes in the mode number space $(l,m)$, all the amplitudes except the modes with $(l,m)=(1,1)$ and $(0,2)$ will be displayed in white, as only these three modes vary with time $t$ in a chaotic manner. 
On the other hand, as time $t$ progresses in the DNS, the modal nonlinear interactions distribute energies to higher-order harmonics and they allow the growth of high-order streamfunction modes such as $\hat{\psi}_{31}$, $\hat{\psi}_{13}$, $\hat{\psi}_{22}$, etc., and high-order temperature modes such as $\hat{\theta}_{11}$, $\hat{\theta}_{04}$, $\hat{\theta}_{31}$, etc. 
As the solution reaches the equilibrium, it is found that the largest amplitudes of the DNS solution are still achieved for the streamfunction mode $\hat{\psi}_{11}=-18.68\mathrm{i}$ and the temperature mode $\hat{\theta}_{02}=-0.3157$ (i.e. $X^{(\mathrm{D})}\simeq12.46$ and $Z^{(\mathrm{D})}\simeq29.75$); however, other high-order modes also have comparably large amplitudes.
It is thus expected that the streamfunction $\psi$ and temperature $\theta$ in the physical space $(x,z)$ are represented not only by the dominant modes with $(l,m)=(1,1)$ and $(l,m)=(0,2)$ but also by other high-order modes. 
In Fig.~\ref{fig:case_r30}(c,d), we also note that the amplitudes in the mode space $(l,m)$ become negligible with amplitudes of order less than $O(10^{-4})$ for $l\geq18$ and $m\geq18$. 
This implies that {GELE} requires the system dimensions with at least $L\simeq18$ and $M\simeq18$ to reproduce the DNS-like results with quantitatively and qualitatively similar nonlinear interactions amongst the high-order modes.

\begin{figure}
\includegraphics[height=7cm]{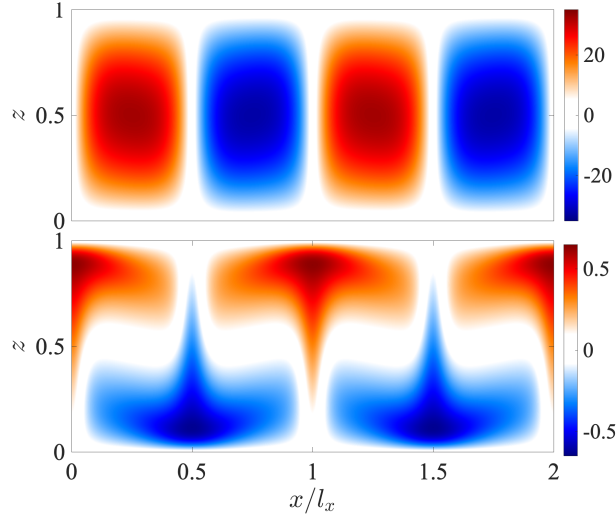}
\caption{\label{fig:psitheta_r30} 
DNS solution of $\psi(x,z)$ (top) and $\theta(x,z)$ (bottom) at the {steady-state} equilibrium at $t=5$ for parameters in Fig.~\ref{fig:case_r30}.
}
\end{figure}
Figure~\ref{fig:psitheta_r30} displays the DNS solution {at the steady-state equilibrium state at $t=5$} in the physical space $(x,z)$ over two streamwise wavelengths (i.e., $x/l_{x}\in [0,2]$).
The streamfunction $\psi$ at the equilibrium represents a pair of vortices (red region: clockwise rotating vortex, blue region: anti-clockwise vortex).
More interestingly, the temperature perturbation $\theta$ exhibits mushroom-shaped convection. 
For both $\psi$ and $\theta$, we see that the dominant spatial periodicity in the $x$-direction is unity.
On the other hand, we see that $\psi(x,z)$ features the spatial periodicity of unity in the $z$-direction while $\theta(x,z)$ shows the spatial periodicity of unity or two depending on the $x$ coordinate. 
These features are captured in the spectral amplitude distributions in Fig.~\ref{fig:case_r30}(c,d) as the most dominant mode in the streamfunction is $\hat{\psi}_{11}$ while both modes $\hat{\theta}_{11}$ and $\hat{\theta}_{02}$ are the most dominant ones for temperature perturbation.
Moreover, the high-order modes also have large amplitudes as we can see a structure like a pointy stem part of the mushroom in the DNS temperature solution $\theta(x,z)$.  
 
\begin{figure*}
\includegraphics[height=6.5cm]{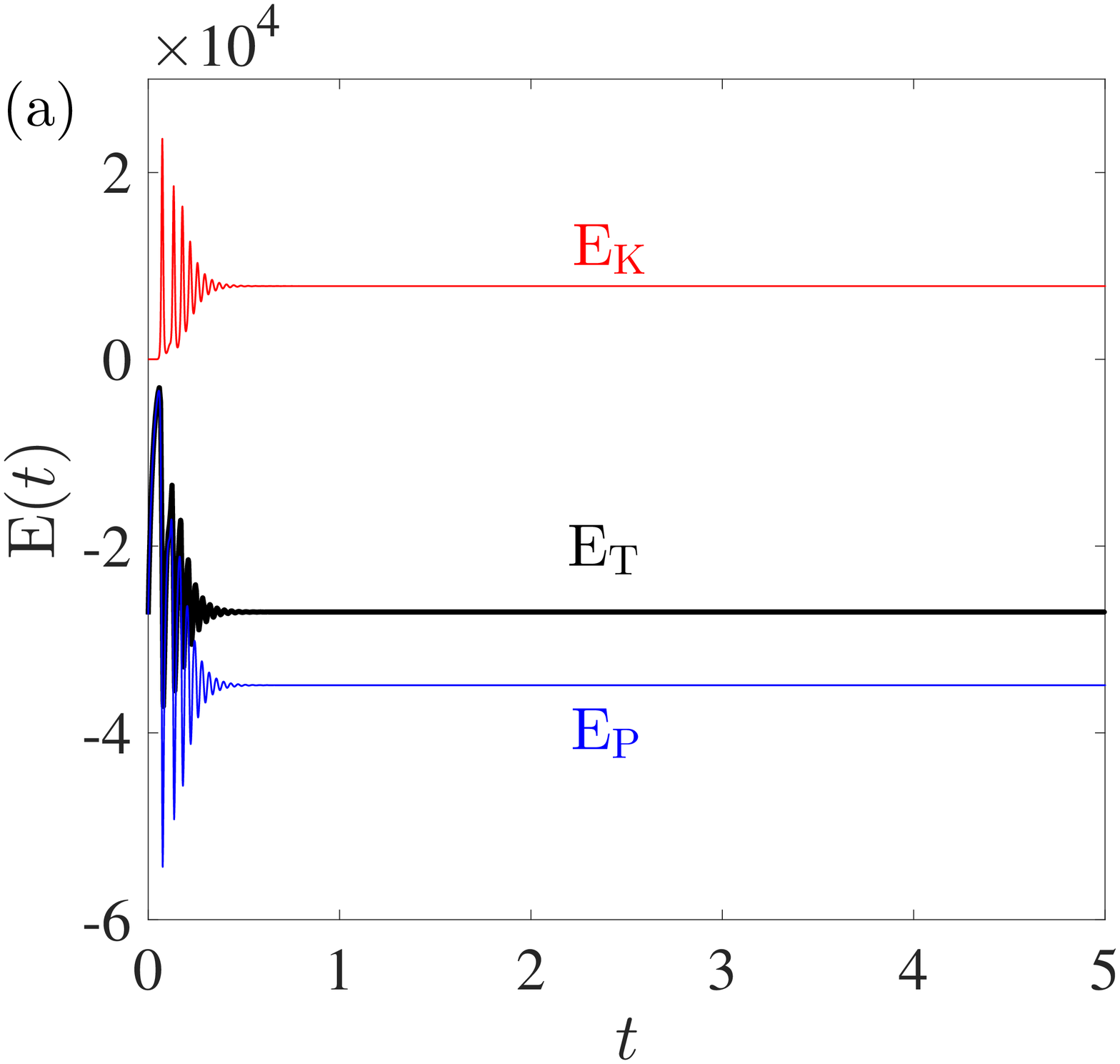}
\includegraphics[height=6.5cm]{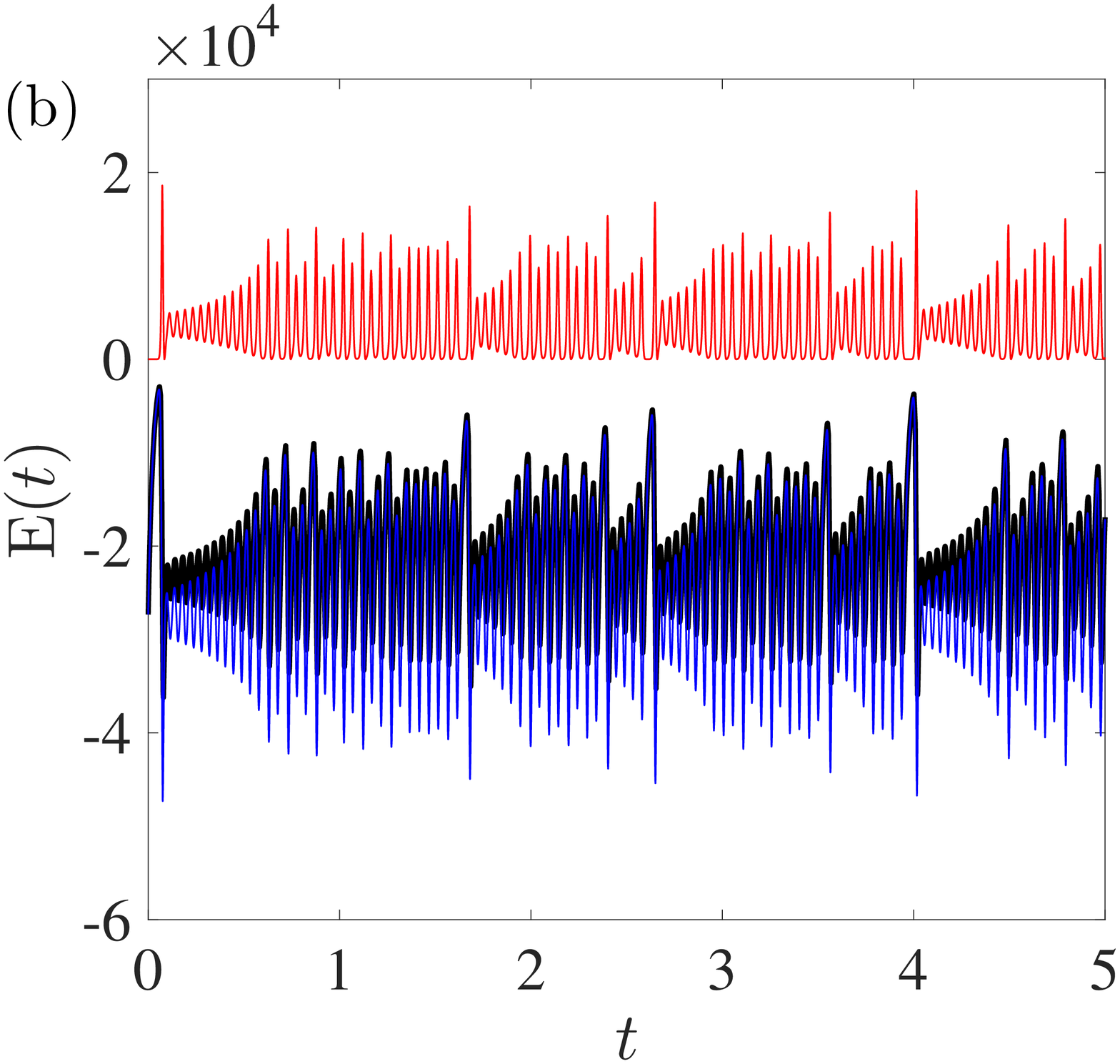}
\includegraphics[height=6.5cm]{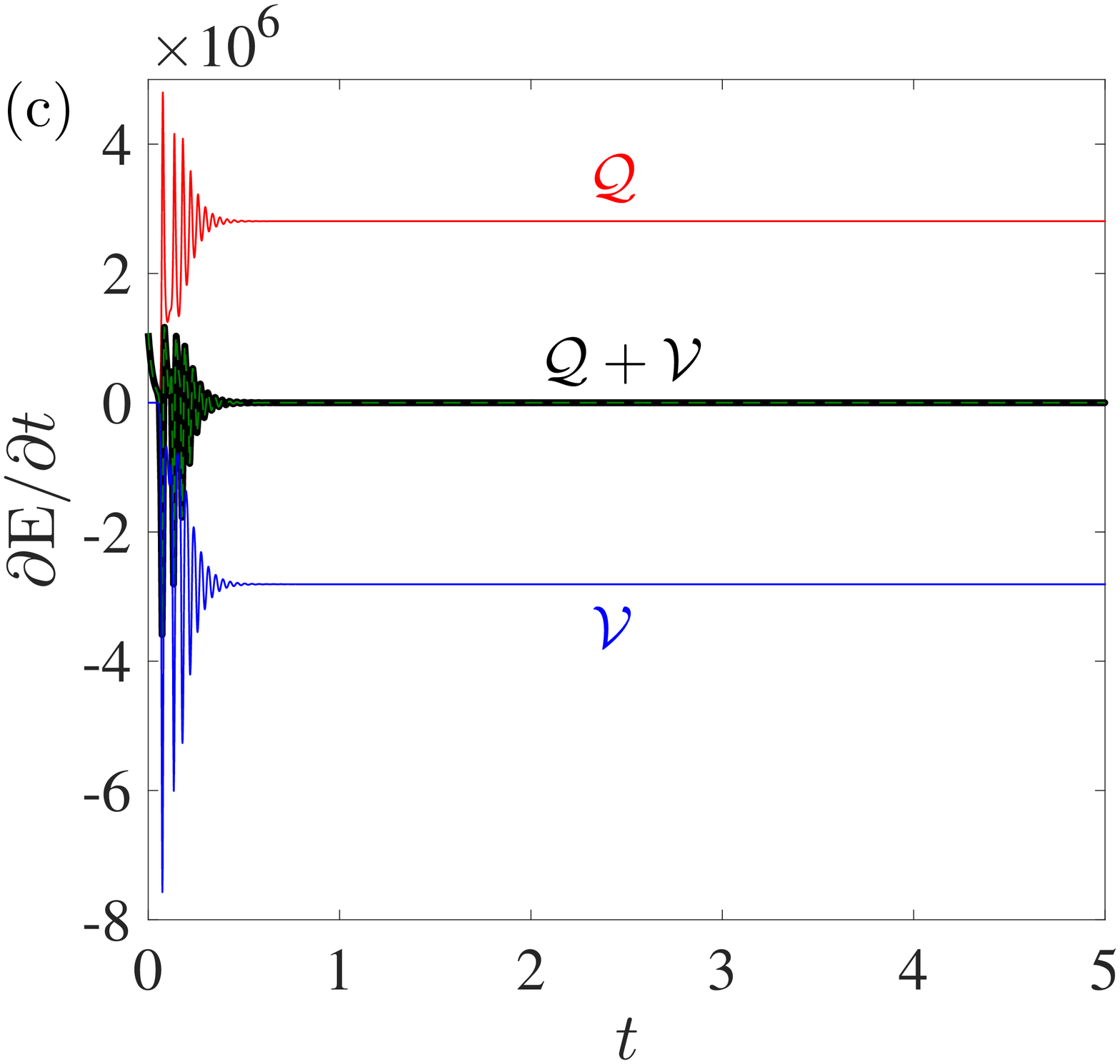}
\includegraphics[height=6.5cm]{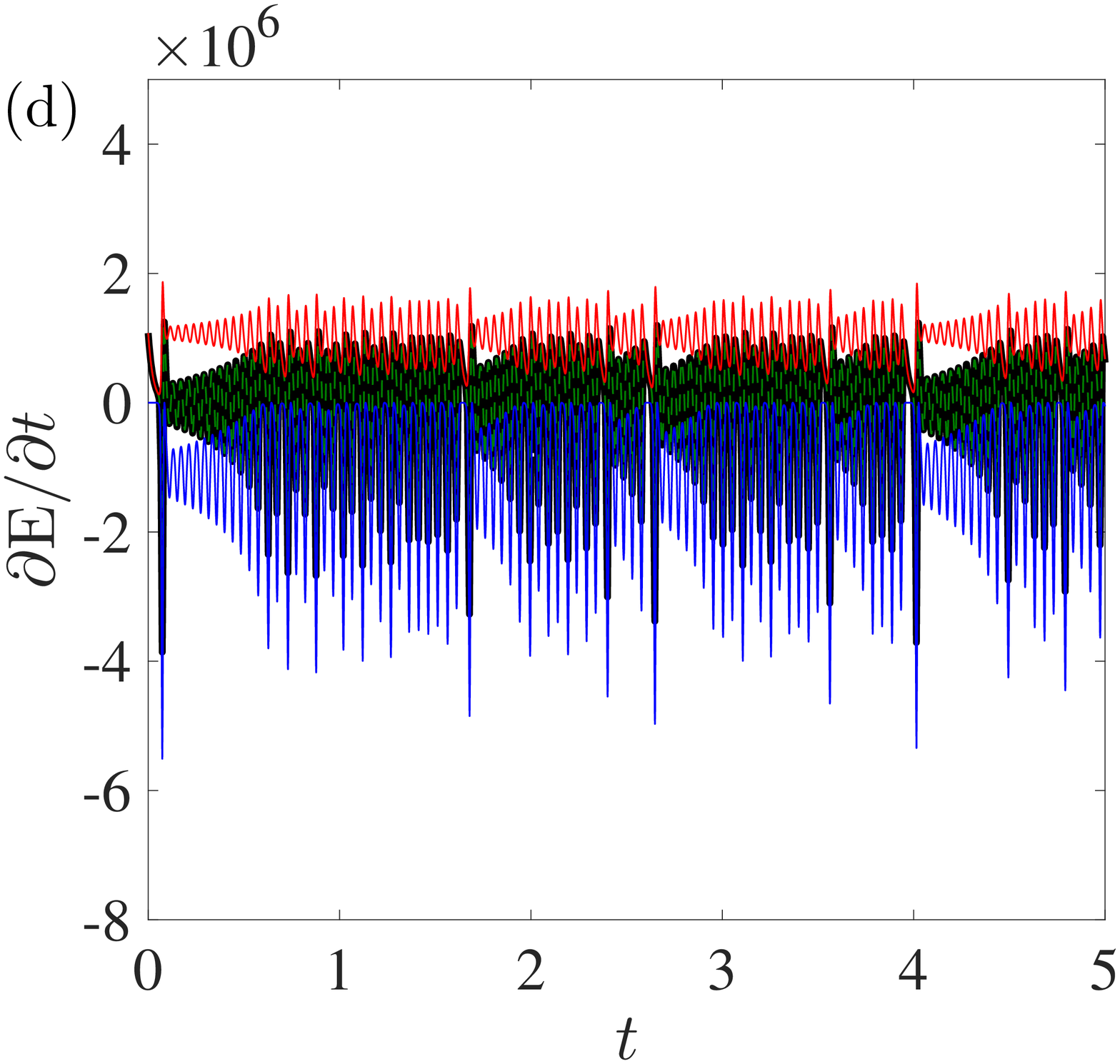}
\caption{\label{fig:case_r30_energy} 
(a,b) Time variation of the total energy $\mathrm{E}_{\mathrm{T}}$ {(black)}, kinetic energy $\mathrm{E}_{\mathrm{K}}$ {(red)}, and potential energy $\mathrm{E}_{\mathrm{P}}$ {(blue)} for the (a) DNS and (b) Lorenz solutions in Fig.~\ref{fig:case_r30}.
(c,d) Time variation of the total energy rate $\partial\mathrm{E}_{\mathrm{T}}/\partial t$ computed directly from $\mathrm{E}_{\mathrm{T}}$ {(green dashed lines overlapped with black solid lines)}, $\mathcal{Q}$ {(red solid lines)}, $\mathcal{V}$ {(blue solid lines)}, and the sum $\mathcal{Q}+\mathcal{V}$ {(black solid lines)} for the (c) DNS and (d) Lorenz solutions. 
}
\end{figure*} 
In Fig.~\ref{fig:case_r30_energy}, we plot the perturbation energy and its time derivative versus time for the DNS and Lorenz solutions of Fig.~\ref{fig:case_r30}.
For both cases, we impose at $t=0$ a small kinetic energy (i.e. $\mathrm{E}_{\mathrm{K}}\simeq4.71\times10^{-3}$) with $X=0.01$. 
And the initial potential energy is negative (i.e. $\mathrm{E}_{\mathrm{P}}\simeq-2.73\times10^{4}$) as the temperature perturbation is stably stratified with $Z=r-1$ at $t=0$.
The total energy $\mathrm{E}_{\mathrm{T}}$ is also negative (i.e. $\mathrm{E}_{\mathrm{T}}\simeq-2.73\times10^{4}$) due to the largely negative potential energy.
Even though the initial kinetic energy is very small, the pair of vortices triggers the instability and the total energy fluctuates with an oscillatory behavior in a transient period, similar to the behavior of $Z(t)$ in Fig.~\ref{fig:case_r30}(a). 
The time variation of the energies for the DNS solution in Fig.~\ref{fig:case_r30_energy}(a) shows the saturation process with the kinetic energy at equilibrium increased from the initial kinetic energy (i.e. the kinetic energy difference $\Delta\mathrm{E}_{\mathrm{K}}=\simeq0.78\times10^{4}$). 
On the other hand, the negative potential energy at the equilibrium is decreased from the initial potential energy (i.e. the potential energy difference $\Delta\mathrm{E}_{\mathrm{P}}\simeq-0.76\times10^{4}$, which implies that the magnitude is increased in the negative direction).
As for the sum, the negative total energy at the equilibrium is slightly increased to $\mathrm{E}_{\mathrm{T}}\simeq-2.71\times10^{4}$ compared to the initial negative total energy (i.e. the increase of the total energy $\Delta\mathrm{E}_{\mathrm{T}}\simeq 2\times10^{2}$, which implies a decrease in magnitude).
The Lorenz solution, on the other hand, does not reach an equilibrium state but it fluctuates in a chaotic manner.
Both the kinetic and potential energies exhibit chaotic temporal variations as shown in Fig.~\ref{fig:case_r30_energy}(b). 
If we average the energies of the Lorenz solution from $t=2$ to $t=5$, we obtain the average total energy $\bar{\mathrm{E}}_{\mathrm{T}}\simeq-2.10\times10^{4}$, the average kinetic energy $\bar{\mathrm{E}}_{\mathrm{K}}\simeq0.32\times10^{4}$, and the average potential energy $\bar{\mathrm{E}}_{\mathrm{P}}\simeq-2.42\times10^{4}$.
While the average kinetic energy of the Lorenz solution is smaller than that of the DNS solution at the equilibrium, the kinetic energy of the Lorenz solution frequently exceeds the equilibrium DNS kinetic energy due to the Lorenz equations' intermittent nature. 

Figure~\ref{fig:case_r30_energy}(c) and (d) display the time derivatives of the energies of the DNS and Lorenz solutions.  
For both solutions, we validate the balance equation (\ref{eq:def_energy_evolution}) by comparing the time derivative $\partial\mathrm{E}_{\mathrm{T}}/\partial t$ directly computed from time-differentiation of $\mathrm{E}_{\mathrm{T}}$ (red dashed line) with the sum $\mathcal{Q}+\mathcal{V}$ (black solid line).
For the DNS solution, the total energy time derivative becomes zero as it reaches the equilibrium and the balance is maintained between the constant negative viscous dissipation $\mathcal{V}$ and the constant positive energy flux $\mathcal{Q}$.
On the other hand, the Lorenz solution does not reach an equilibrium as the viscous dissipation $\mathcal{V}$ and the energy flux $\mathcal{Q}$ do not balance but they fluctuate with time in a chaotic manner; therefore, the time derivative of the total energy $\partial\mathrm{E}_{\mathrm{T}}/\partial t$ for the Lorenz solution never stays at zero.

\subsection{Connection between Lorenz and DNS solutions}
\begin{figure}
\includegraphics[height=5.8cm]{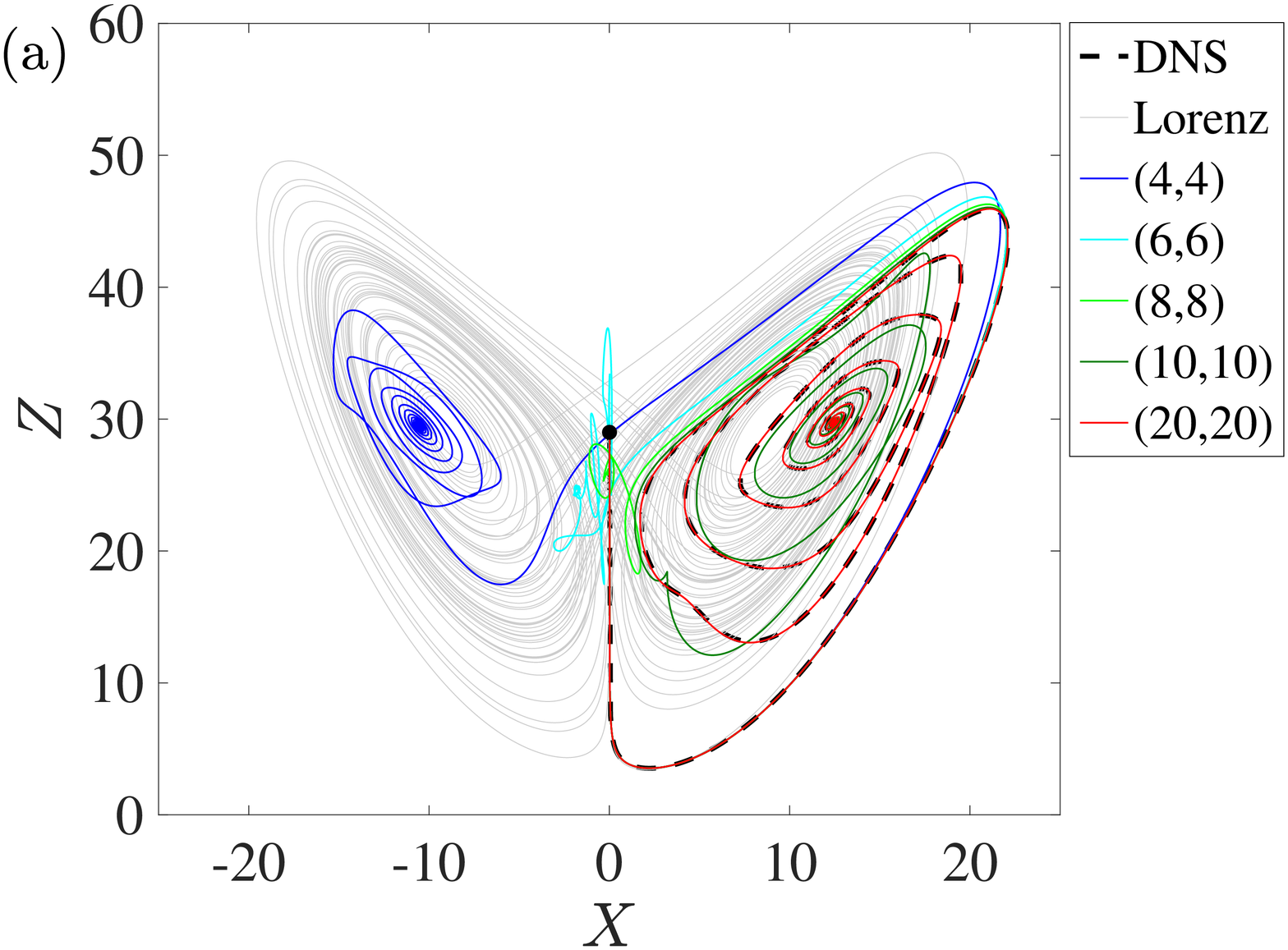}
\includegraphics[height=5.8cm]{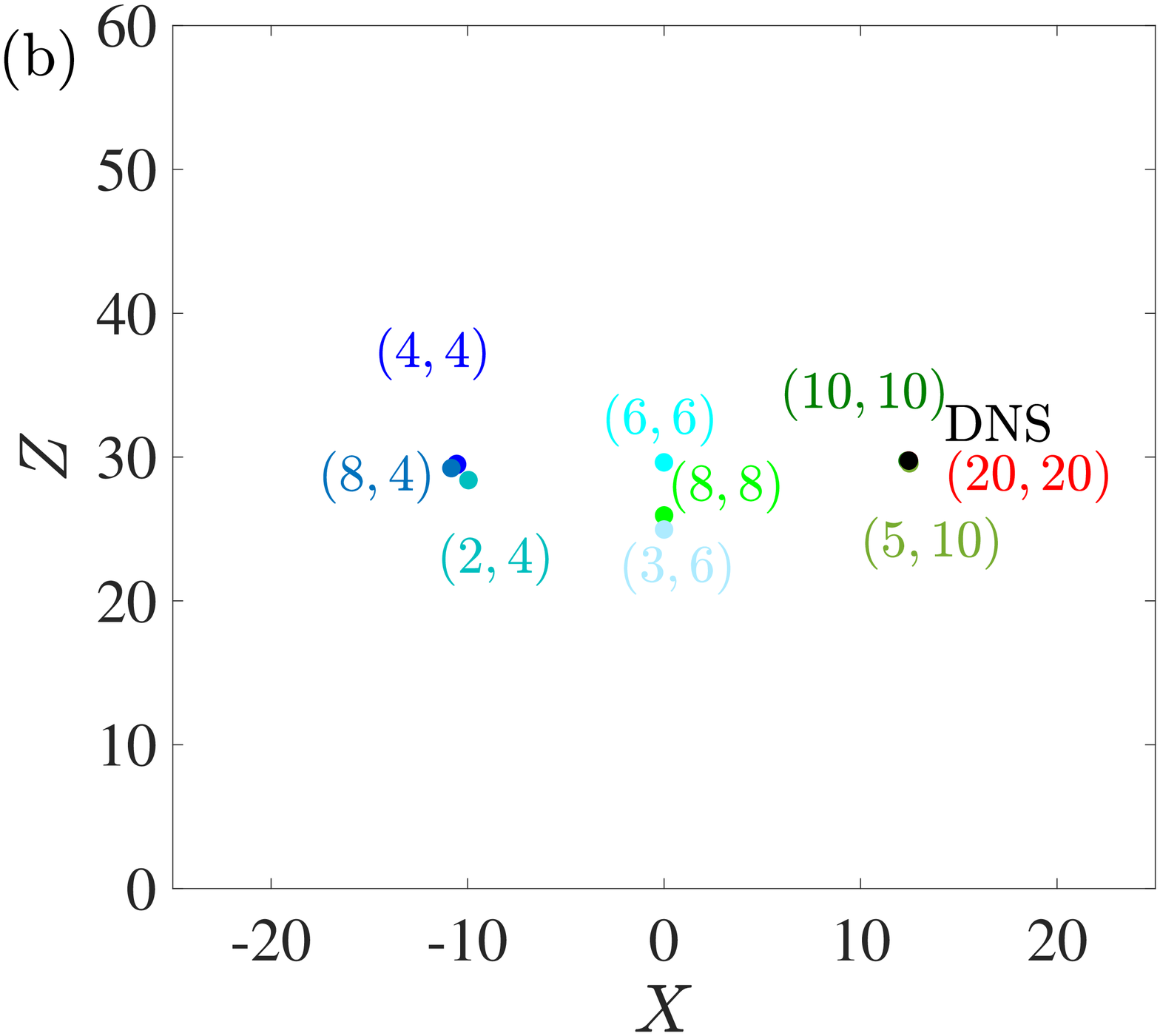}
\caption{\label{fig:XZ_r30_GNDS} 
{(a)} Trajectories on the $(X,Z)$-plane for various {solutions of the GELE} with different $L$ and $M$ (colored solid lines), the Lorenz equations (gray solid line), and the DNS solution (black dashed line) at $r=30$. Black circle indicates the initial condition $(X,Z)=(0.01,29)$.
{(b) Various fixed points for converging solutions of the GELE with different $(L,M)$ and DNS solution in (a).}
}
\end{figure} 
In this subsection, we now investigate with {GELE} how solutions transition from the Lorenz equations to the DNS as the mode limits $L$ and $M$ are increased. 
Given the same initial condition $(X,Y,Z)=(0.01,0,r-1)$, Fig.~\ref{fig:XZ_r30_GNDS}{(a)} shows trajectories on the $(X,Z)$-plane of solutions with various values of $L$ and $M$.
The trajectories of the DNS and Lorenz solutions are the same as the ones in Fig.~\ref{fig:case_r30}(b), only displayed with different line styles in Fig.~\ref{fig:XZ_r30_GNDS}.
It is remarkable that the high-order solutions other than the Lorenz solution do not exhibit chaotic attractors but converge to fixed points; for instance, the trajectories converge to $(X^{(\mathrm{G})},Z^{(\mathrm{G})})\simeq(-10.55,29.49)$ for $(L,M)=(4,4)$, $(X^{(\mathrm{G})},Z^{(\mathrm{G})})\simeq(-0.006,29.63)$ for $(L,M)=(6,6)$, $(X^{(\mathrm{G})},Z^{(\mathrm{G})})\simeq(0,25.94)$ for $(L,M)=(8,8)$, $(X^{(\mathrm{G})},Z^{(\mathrm{G})})\simeq(12.39,29.75)$ for $(L,M)=(10,10)$, and $(X^{(\mathrm{G})},Z^{(\mathrm{G})})\simeq(12.46,29.75)$ for $(L,M)=(20,20)$.
Fixed points of the {GELE solutions} depend on $L$ and $M$ {as shown in Fig.~\ref{fig:XZ_r30_GNDS}(b)}, but it is verified that they approach the fixed points of the DNS as $L$ and $M$ increase. 
The trajectory of the system with $(L,M)=(10,10)$ is slightly different from the trajectory of the DNS solution in the transient period, but the final fixed point $(X^{(\mathrm{G})},Z^{(\mathrm{G})})\simeq(12.39,29.75)$ is very similar to the equilibrium $(X^{(\mathrm{D})},Z^{(\mathrm{D})})\simeq(12.46,29.75)$ of the DNS solution. 
For higher orders of $L>10$ and $M>10$, the trajectories of the {GELE} solution become equivalent to those of the DNS solution. 
As the system order increases, the number of possible fixed points increases and onto which fixed point a trajectory settles depends on the initial condition. 
We have checked that the same initial condition for different $L$ and $M$ leads to the same fixed point when $L$ and $M$ are sufficiently large. 
Further discussion on the initial-condition dependency will be provided in another subsection.

\begin{figure}
\includegraphics[height=10cm]{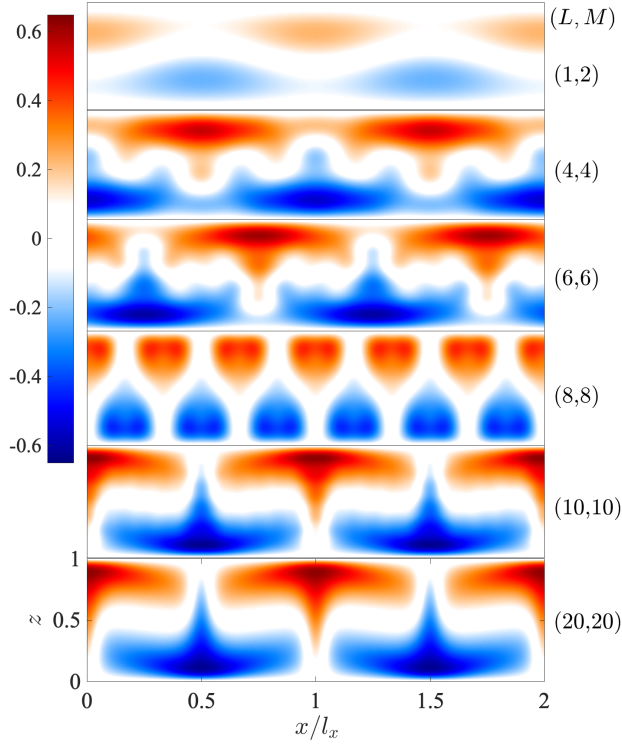}
\caption{\label{fig:theta_r30} 
Temperature perturbation $\theta(x,z)$ at $t=5$ obtained from {GELE} for various sets of $(L,M)$ and parameters in Fig.~\ref{fig:case_r30}.
}
\end{figure}
To understand in a more visual way how a solution transitions from the Lorenz equations to the DNS, Fig.~\ref{fig:theta_r30} shows temperature perturbation $\theta(x,z)$ over two streamwise wavelengths $2l_{x}$ for the {GELE} solutions with various sets of $(L,M)$.
Only the Lorenz solution with $(L,M)=(1,2)$ at the top of Fig.~\ref{fig:theta_r30} is not at equilibrium at $t=5$ as the Lorenz solution lies on a chaotic attractor before and after $t=5$, while other {GELE} solutions of higher orders reach their equilibrium states. 
For all solutions in Fig.~\ref{fig:theta_r30}, we recognize that the dominant spatial periodicity in the $z$-direction is two (i.e. the dominant mode number is $m=2$).
On the other hand, the dominant spatial periodicity in the $x$-direction varies with the system orders $L$ and $M$.
For instance, the temperature perturbations for $(L,M)=(4,4)$ and $(6,6)$ show a wiggly pattern around the center line $z=0.5$ and it is difficult to determine by inspection which mode number $l$ is the dominant one. 
For the temperature perturbation of $(L,M)=(8,8)$, it is noticeable that the dominant periodicity in the $x$-direction is $l=3$ {(i.e. the dominant wavelength is $l_{x}/3$). 
A similar structure with the dominant spatial periodicity $l=3$ is observed for the case $(L,M)=(3,6)$ (not shown) when the same initial condition is imposed.}
As the system limits $L$ and $M$ are further increased, the {GELE} equilibrium solutions for $L\geq10$ and $M\geq10$ become equivalent to the DNS solution in Fig.~\ref{fig:psitheta_r30}.

\subsection{Periodic and chaotic solutions}
\begin{figure}
\includegraphics[height=7cm]{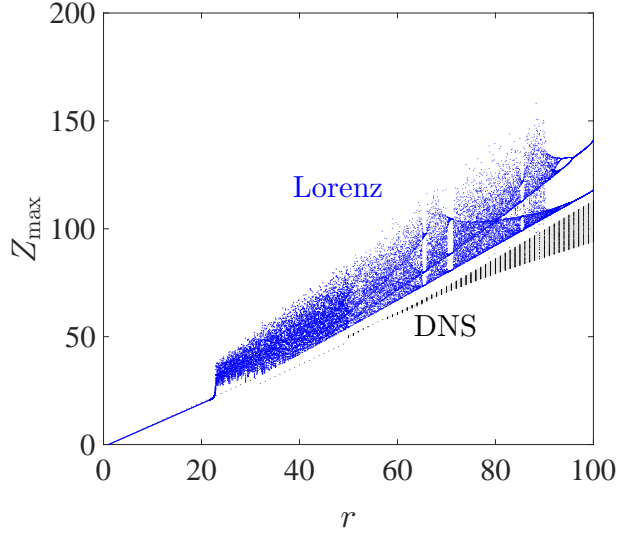}
\caption{\label{fig:Zmax_Lorenz_DNS} 
Bifurcation diagrams of $Z_{\max}$ versus $r$ for the Lorenz (blue) and DNS (black) solutions. 
Dots denote actual $Z_{\max}$ picked up at each local maximum, and gray area denotes the possible range of $Z_{\max}$ due to the appearance of the limit tori for $r\geq58$.
{For the DNS solutions, the interval $\Delta r=1$ is used.}
}
\end{figure}
We now investigate how the solution behaviors change as $r$ is increased.
For each $r$, we still use the Lorenz-like initial condition with $(X,Y,Z)=(0.01,0,r-1)$ and other variables set to zero. 
In Fig.~\ref{fig:Zmax_Lorenz_DNS}, we plot the bifurcation diagrams of $Z_{\max}$ versus $r$ for the Lorenz and DNS solutions. 
The local maxima of $Z$, $Z_{\max}$, are picked up after truncation of the transient period ($0\leq t\leq 3$) from the solution\citep[][]{Yu1996,Park2015PS}, and we define hereafter the $Z$-periodicity of the solution as the number of $Z_{\max}$. 
{Integer choices in $r$ with the} interval $\Delta r=1$ is used to plot the bifurcation diagram of the DNS solution.
Our focus is not on the blue-dotted Lorenz bifurcation, which has already been investigated extensively in previous studies (see e.g. \citet{Dullin2007}), but on the bifurcation behavior of the DNS solution in the parameter space $r$. 
While the Lorenz equations bifurcate beyond $r>24$, the trajectories of DNS solutions converge to fixed points in the range $1<r<50$. 
The DNS bifurcation curve is slightly dropped in the range $30<r<50$ due to the convergence to a fixed solution of the streamwise periodicity of $3$ in this particular range of $r$, while the solutions in the range $r\leq 30$ have the streamwise periodicity of unity as shown in Fig.~\ref{fig:psitheta_r30} for $r=30$. 
Beyond $r\geq 50$, it is found that limit cycles with the $Z$-periodicity of unity appear in the range $50\leq r\leq 58$ and limit tori appear for $r>58$.
For a limit torus, it is thought that there are infinitely many distinct $Z_{\max}$, so we have the gray shaded area in Fig.~\ref{fig:Zmax_Lorenz_DNS} indicating the possible range of $Z_{\max}$. 
We see that the width of the gray area increases gradually as $r$ increases. 

\begin{figure*}
\includegraphics[height=5.5cm]{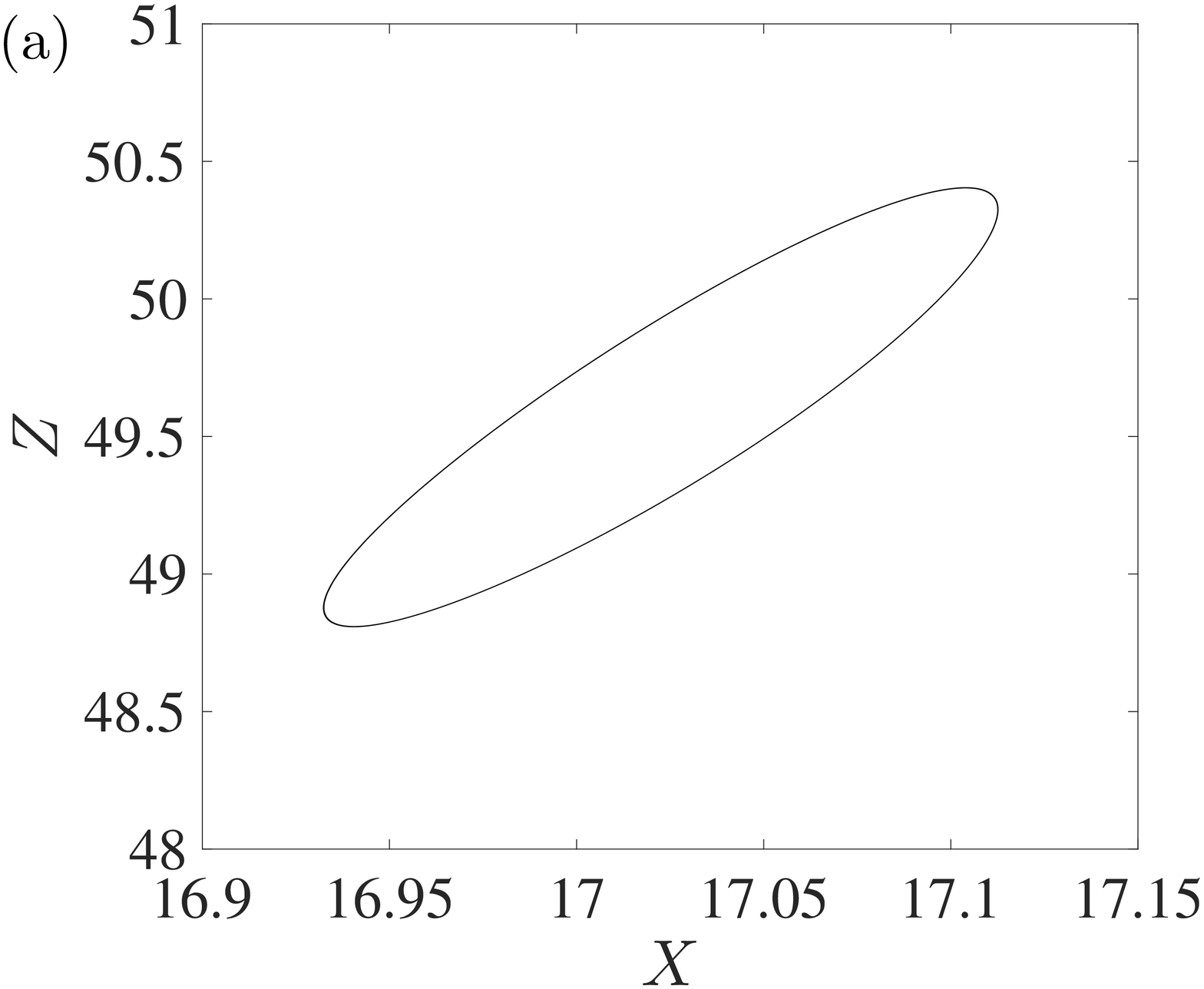}
\includegraphics[height=5.5cm]{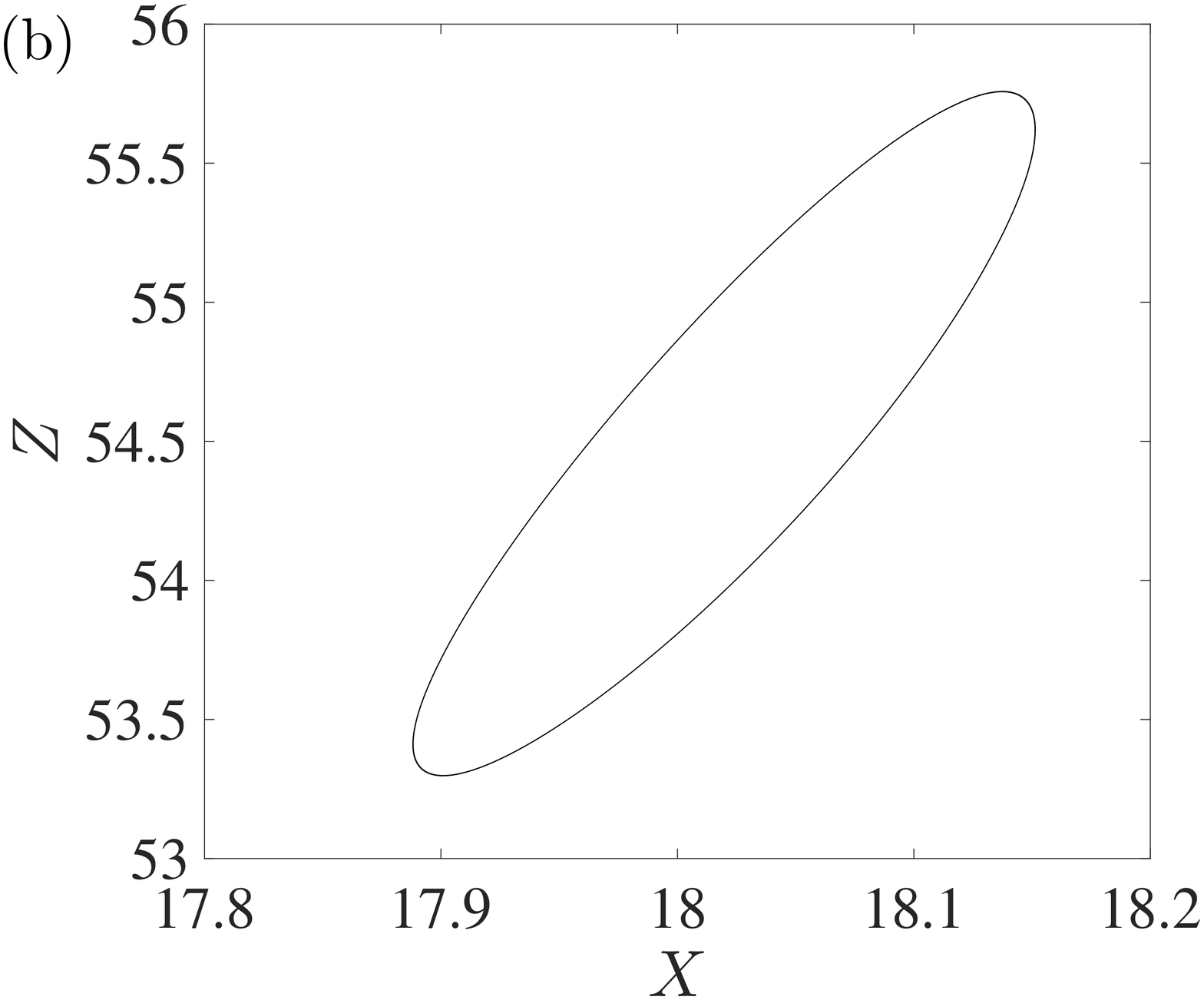}
\includegraphics[height=5.5cm]{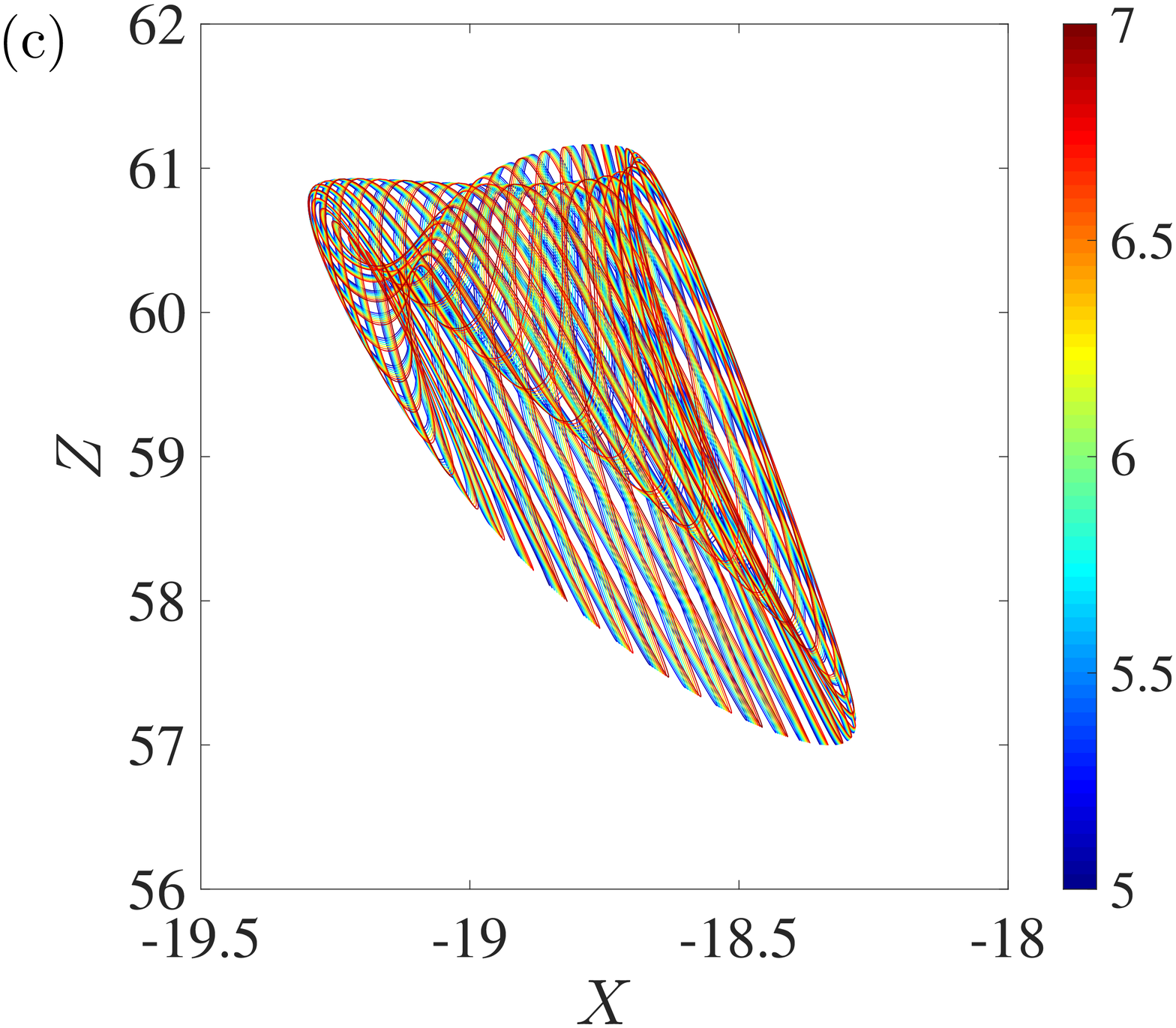}
\includegraphics[height=5.5cm]{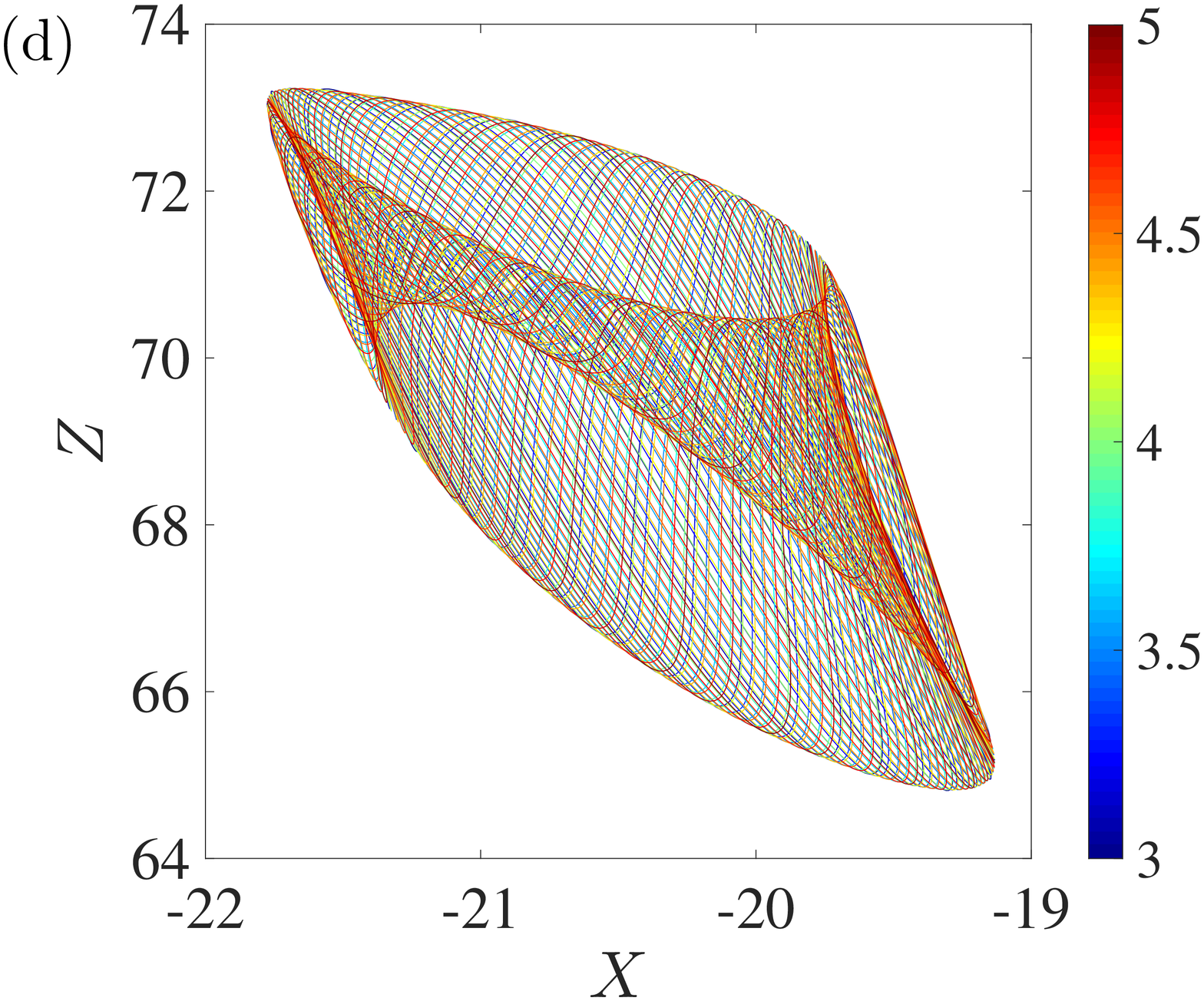}
\caption{\label{fig:cases_XZ} 
Trajectories on the $(X,Z)$-plane computed from the DNS for (a) $r=50$, (b) $r=55$, (c) $r=60$, and (d) $r=70$. 
{In (c) and (d), the changing colors of the limit tori are based on time $t$ as displayed in the colorbars.}
}
\end{figure*} 
To see more clearly what types of periodic solutions are observed, we show in Fig.~\ref{fig:cases_XZ} the trajectories of the DNS solutions on the $(X,Z)$-plane.
In the range $1<r<50$, it is verified that the DNS solution saturates nonlinearly and its trajectory converges to a fixed solution as reaching the equilibrium state.
If we plot only the fixed solution on the $(X,Z)$-plane, it will appear as a dot.
As $r$ increases further, in the range $50\leq r\leq 58$, the DNS solution becomes periodic and the solution exhibits a limit cycle with the $Z$-periodicity of unity as shown in Fig.~\ref{fig:cases_XZ}(a,b) for $r=50$ and 55. 
As $r$ increases beyond $r=58$, the solution's trajectory no longer lies on a limit cycle; for instance, the trajectory in Fig.~\ref{fig:cases_XZ}(c) at $r=60$ does not exhibit a limit cycle of the $Z$-periodicity of unity on the $(X,Z)$-plane.
The trajectory is, however, somehow regular and bounded.
A more regular pattern is observed for the trajectory at $r=70$ as shown in Fig.~\ref{fig:cases_XZ}(d). 

\begin{figure*}
\includegraphics[height=5.2cm]{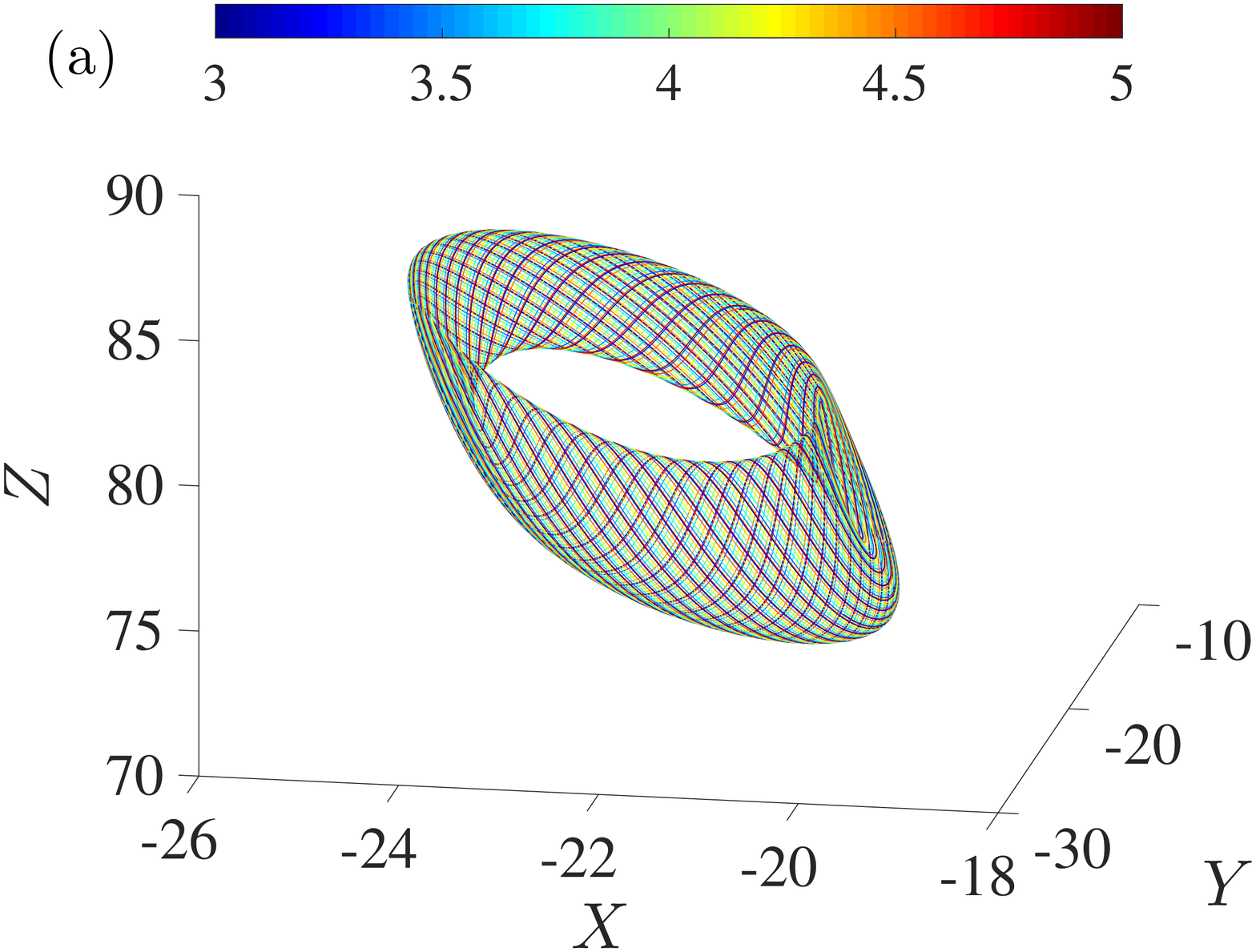}
\includegraphics[height=5.2cm]{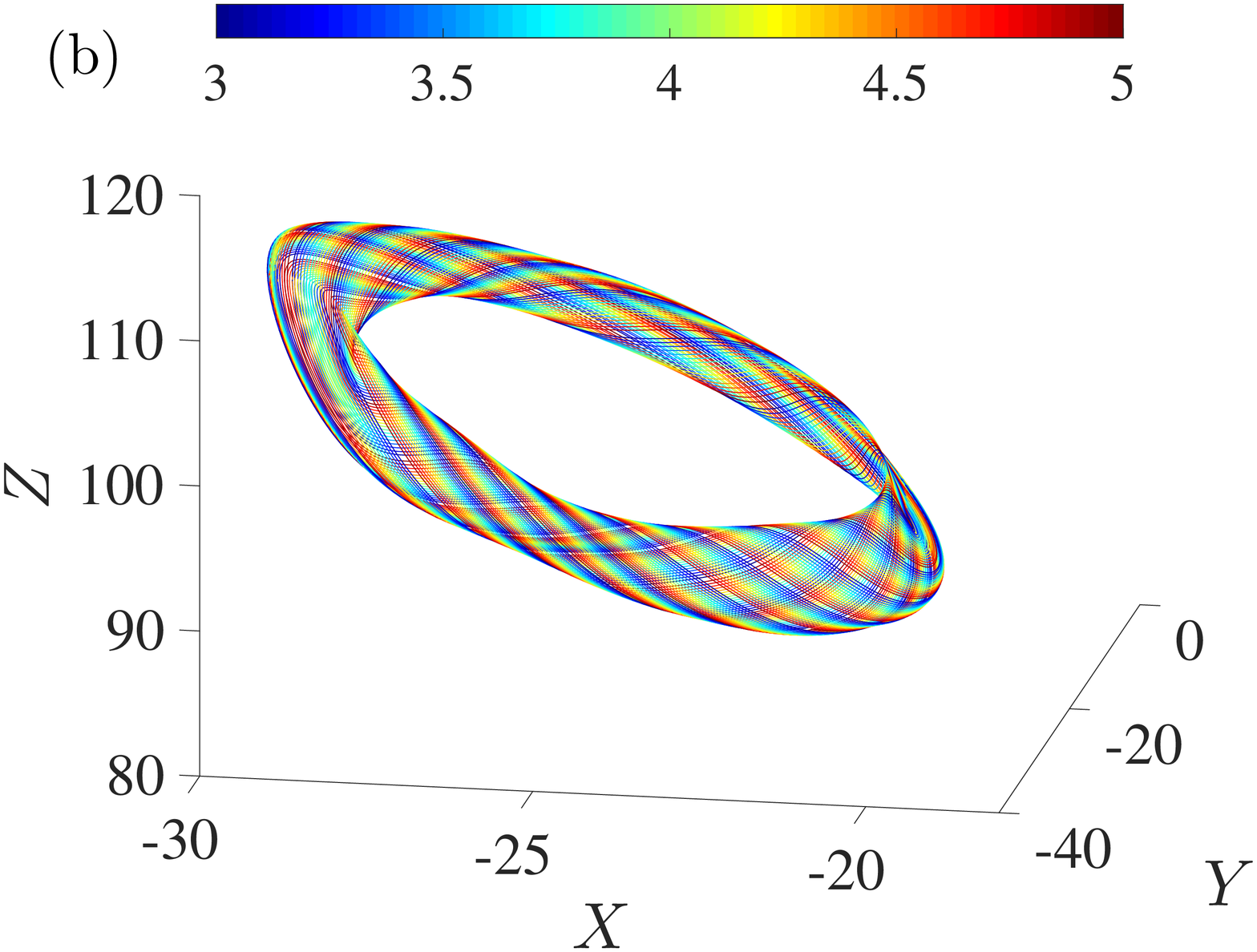}
\includegraphics[height=5.2cm]{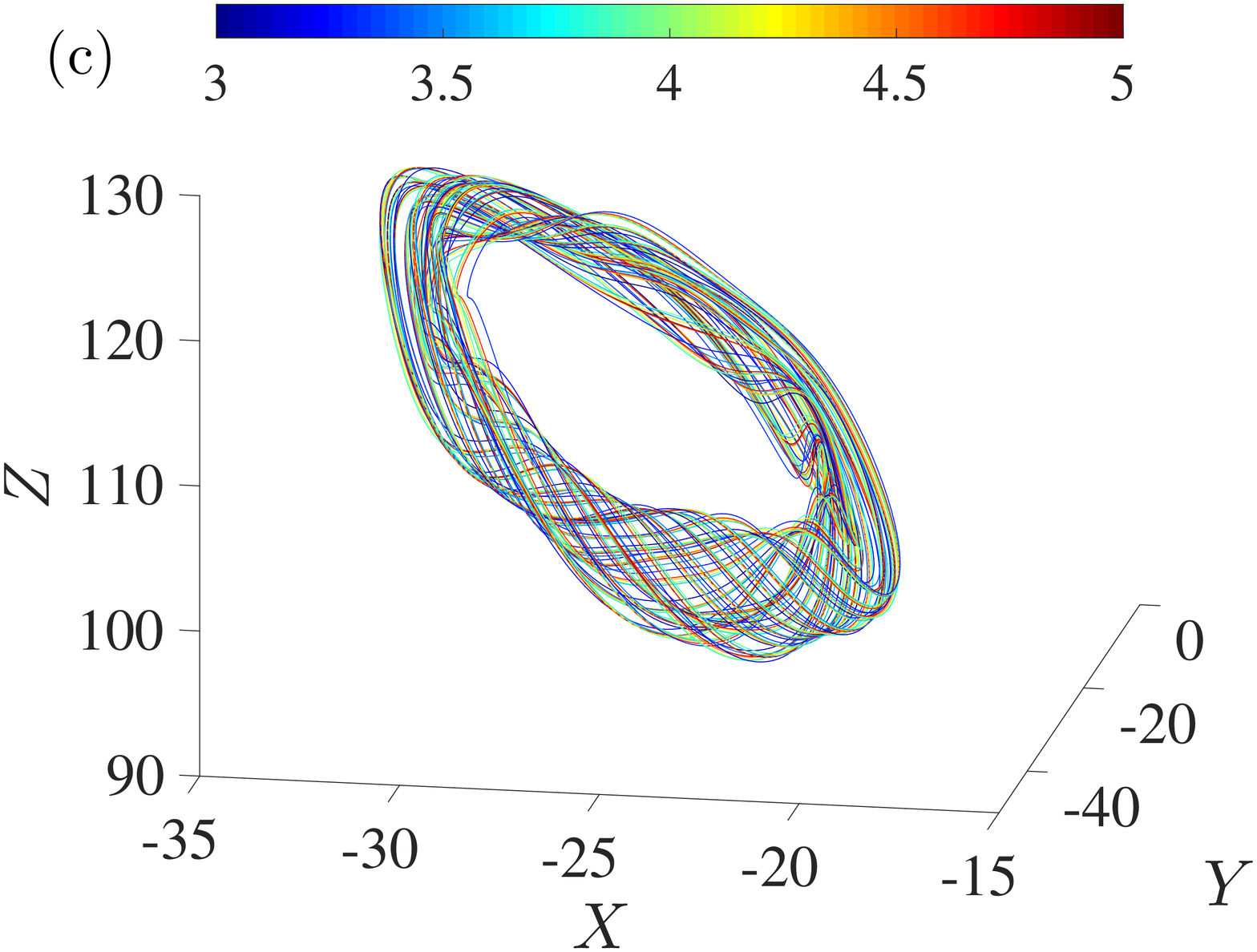}
\includegraphics[height=5.2cm]{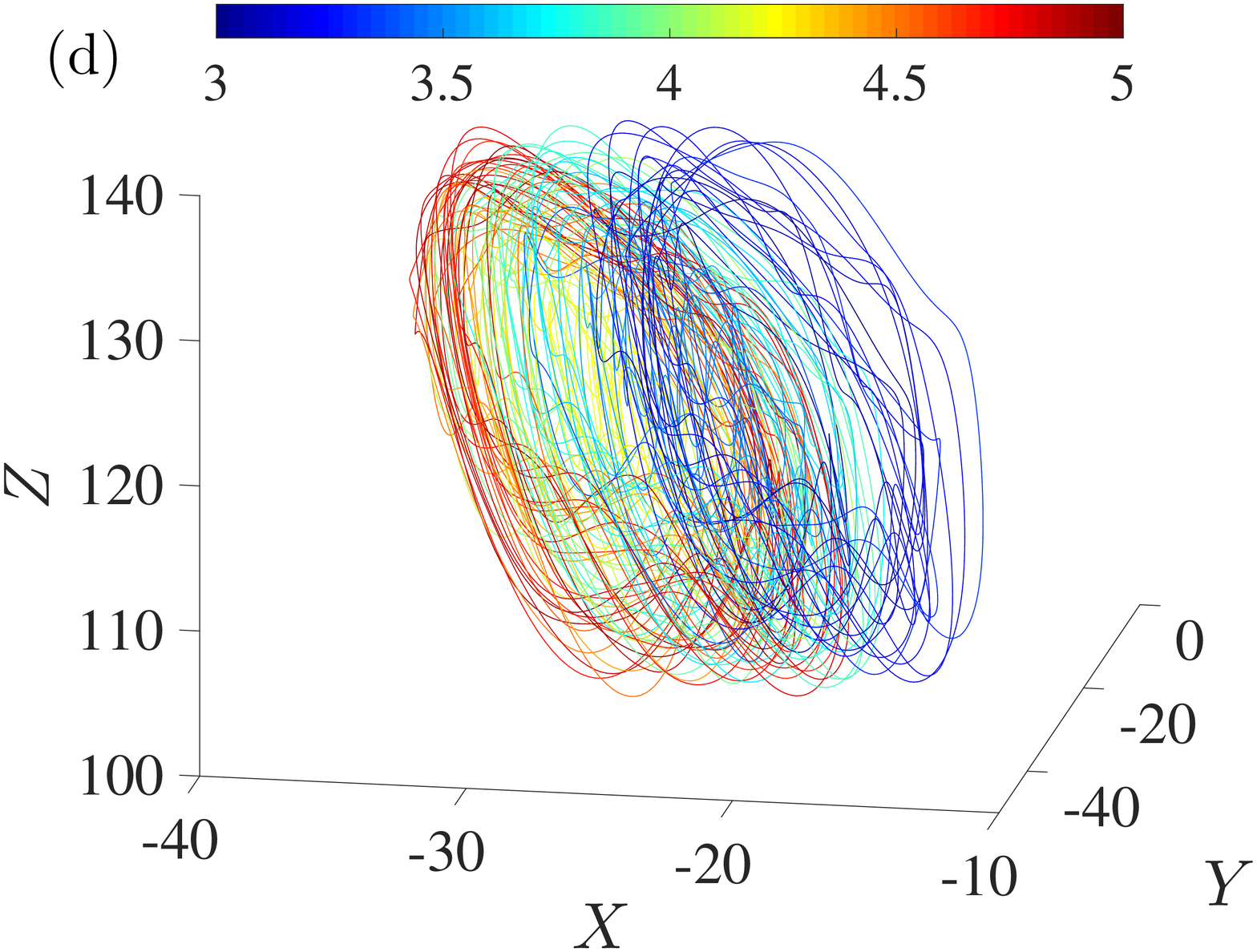}
\includegraphics[height=5.2cm]{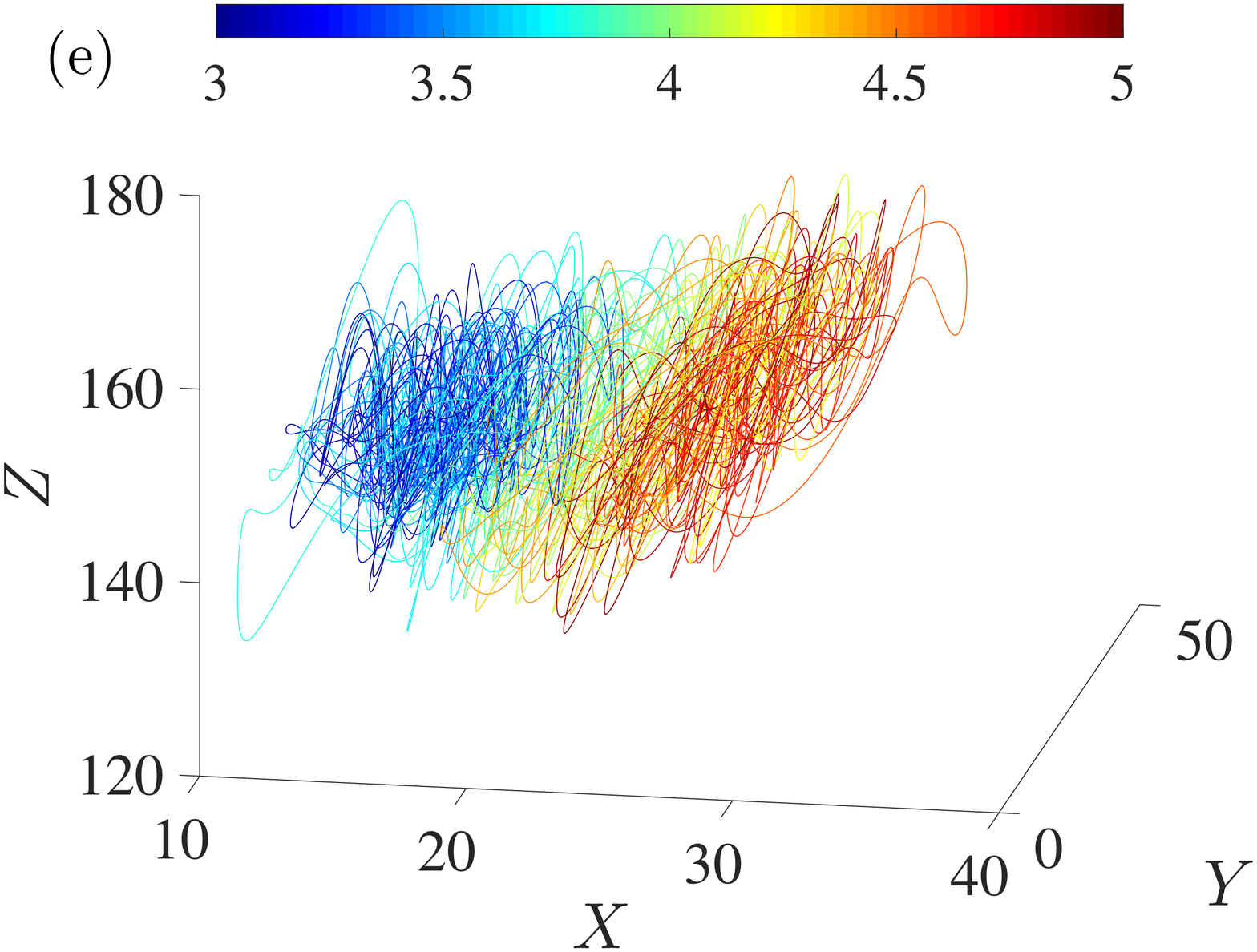}
\includegraphics[height=5.2cm]{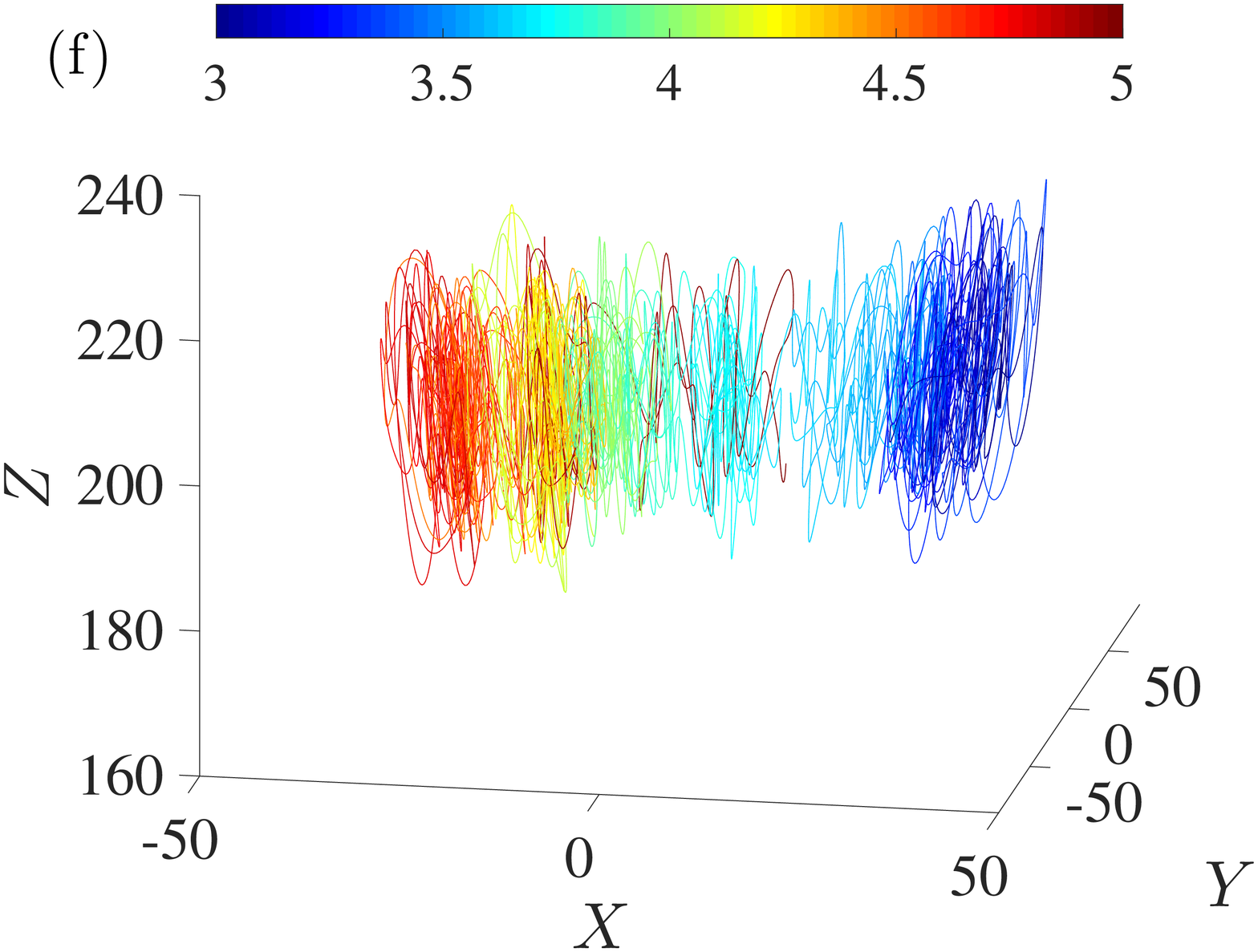}
\caption{\label{fig:cases_XYZ} 
Trajectories of the DNS solutions in the $(X,Y,Z)$-space for (a) $r=80$, (b) $r=100$, (c) $r=110$, (d) $r=120$, (e) $r=150$, and (f) $r=200$. 
{Colorbars display the value of time $t$ corresponding to each color of the trajectories.}
}
\end{figure*} 
To better understand the bounded trajectories in the range $r>58$, we plot in Fig.~\ref{fig:cases_XYZ} three-dimensional trajectories of the DNS solutions in the $(X,Y,Z)$-space for various values of $r$ where the solution no longer lies on a limit cycle and does not converge to a fixed point.
At $r=80$ as shown in Fig.~\ref{fig:cases_XYZ}(a), the solution lies on a smooth limit torus, which is known to be observed in the presence of quasiperiodicity \citep[][]{Grebogi1985}. 
It is verified that trajectories of the solutions in the range $58< r<80$ (including the ones at $r=60$ and $r=70$ shown in Fig.~\ref{fig:cases_XZ}c and d) also lie on limit tori. 
The solution at $r=100$ in Fig.~\ref{fig:cases_XYZ}(b) exhibits a limit torus attractor as well, but it is now twisted along the toroidal direction. 
The solution's irregularity becomes more apparent as $r$ increases further.
At $r=110$, the trajectory has an irregular torus shape (Fig.~\ref{fig:cases_XYZ}(c)), that is, the solution does not exhibit any regular-shape attractor (e.g. limit cycles, limit tori). The trajectory continues to move irregularly as $r\geq120$ (see Fig.~\ref{fig:cases_XYZ}(d)--(f)).
It is noticeable that such irregular chaotic solutions cover wider ranges of $(X,Y,Z)$ in the phase space as $r$ increases. 

\begin{figure*}
\includegraphics[height=5.7cm]{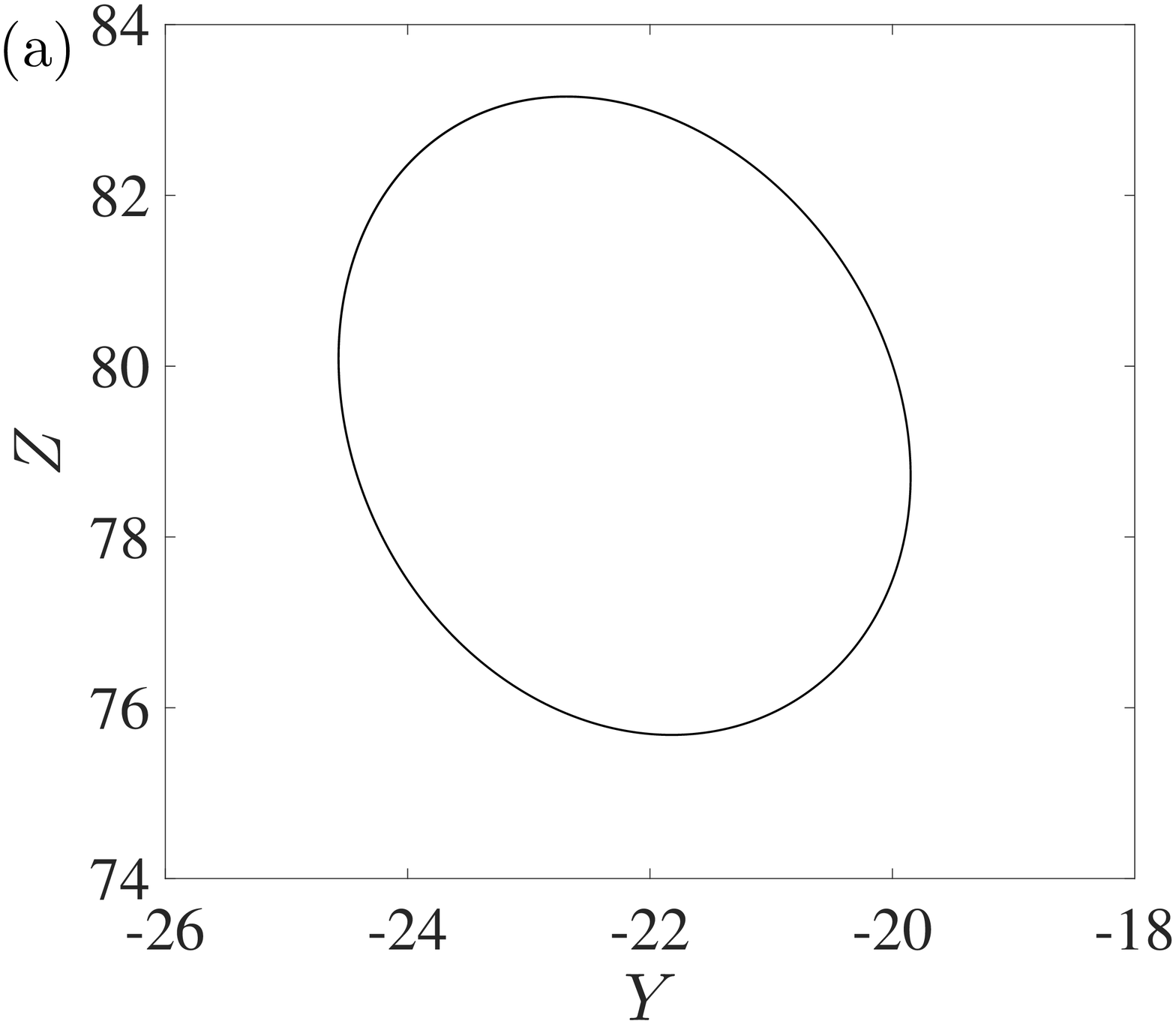}
\includegraphics[height=5.7cm]{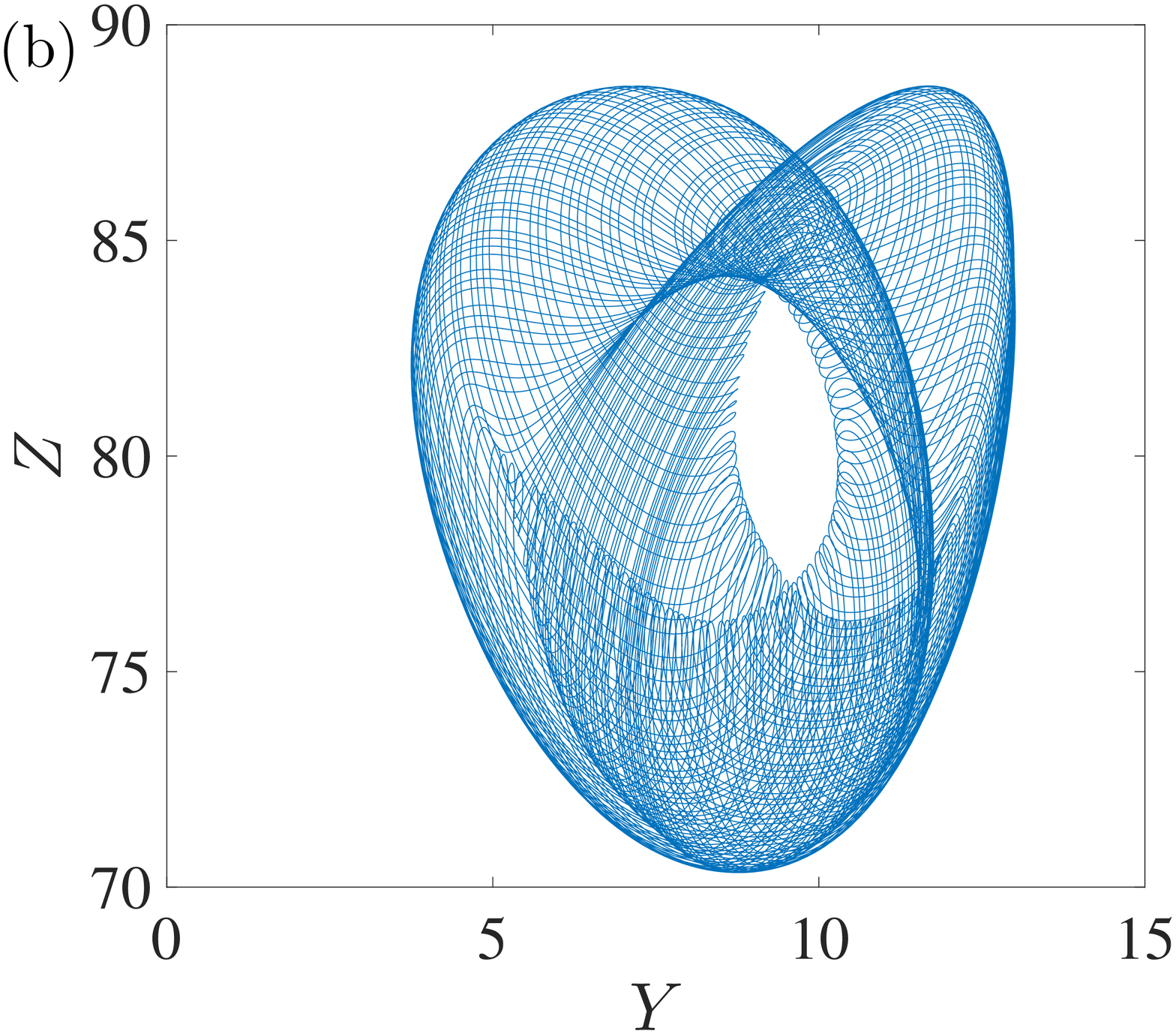}
\includegraphics[height=5.9cm]{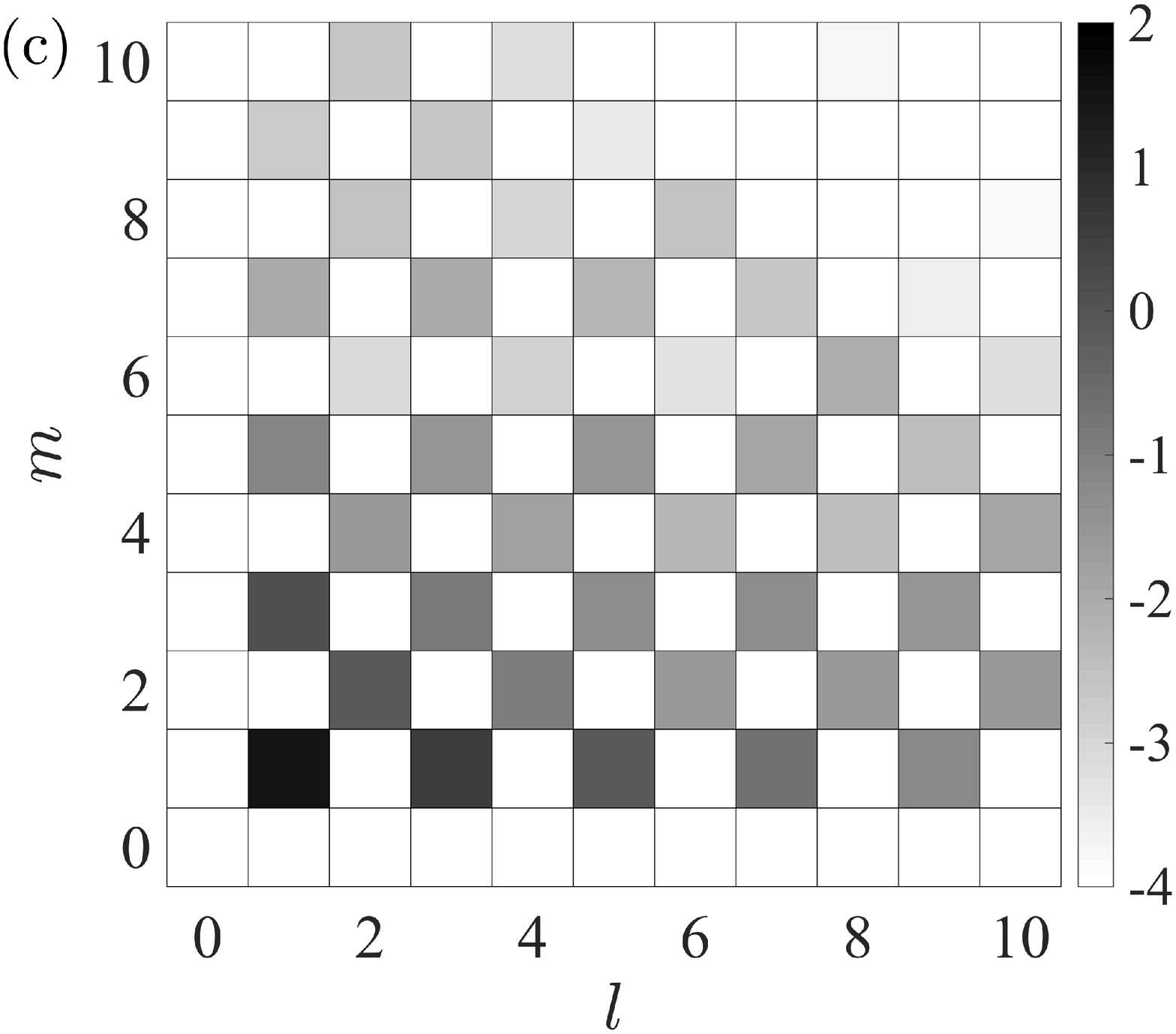}
\includegraphics[height=5.9cm]{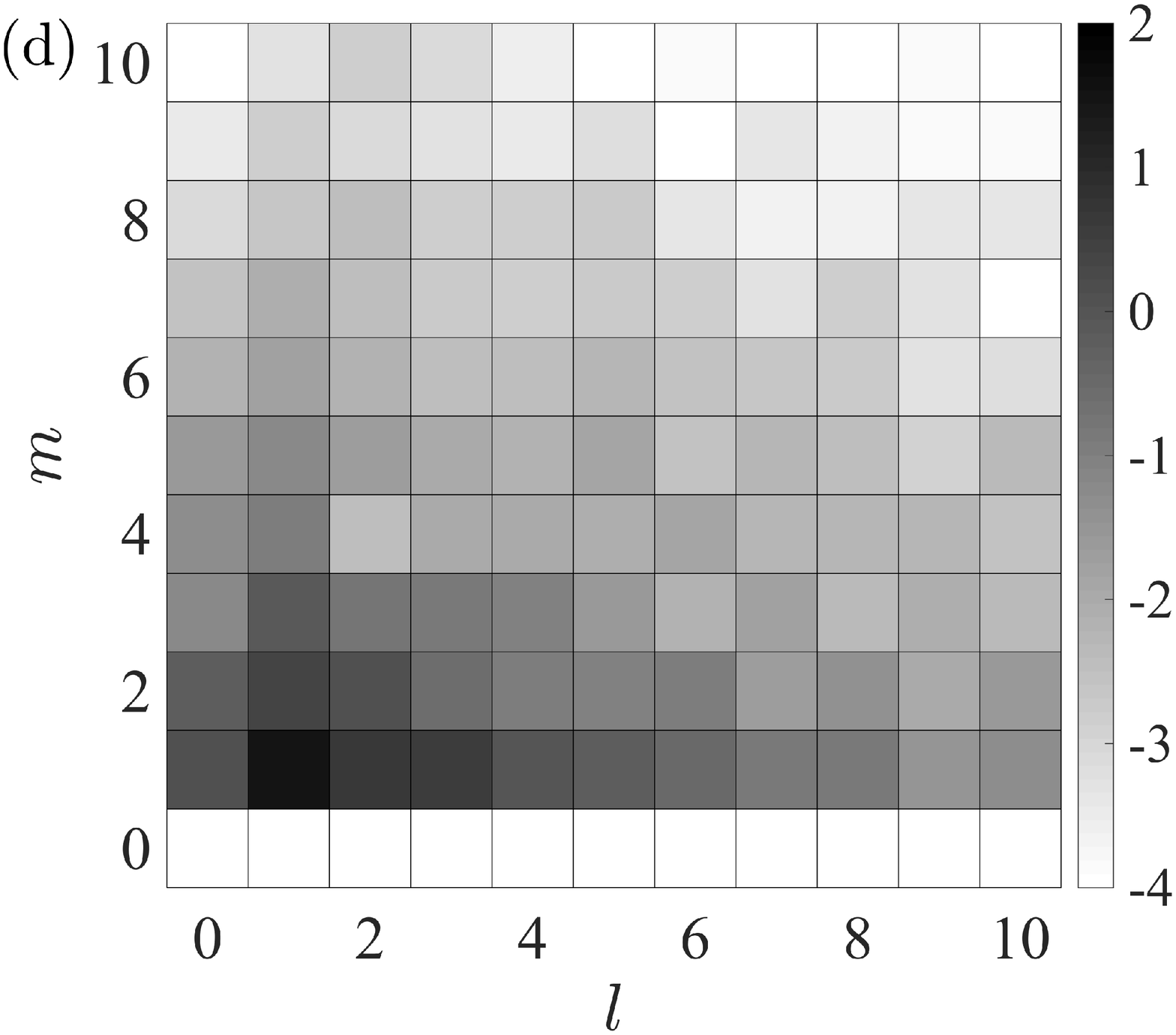}
\caption{\label{fig:cases_L10M10YZ} 
(a,b) Trajectories of the {GELE} solutions on the $(Y,Z)$-plane after a transient time period for $r=80$, $(L,M)=(10,10)$ from (a) Lorenz-like and (b) random initial conditions.
(c) Distribution of the amplitude $\log_{10}(|\hat{\psi}_{lm}|)$ in the parameter space $(l,m)$ for a {GELE} solution on the black limit cycle in (a). (d) The amplitude distribution $\log_{10}(|\hat{\psi}_{lm}|)$ for a {GELE} solution on the blue limit torus in (b). 
}
\end{figure*} 
To verify if a limit torus is also observable in the {GELE}, we compute the solutions of the {GELE} of orders $(L,M)=(10,10)$ at $r=80$ (Fig.~\ref{fig:cases_L10M10YZ}).
It is found that, if the Lorenz-like initial condition (i.e. $(X,Y,Z)=(0.01,0,r-1)$ and other variables zero) is imposed, the {GELE} solution lies on a limit cycle as shown in Fig.~\ref{fig:cases_L10M10YZ}(a), which is different from the DNS solution's limit torus behavior.
To understand this different outcome, we plot the amplitude $\hat{\psi}_{lm}$ in the parameter space $(l,m)$ in Fig.~\ref{fig:cases_L10M10YZ}(b), and we see that the limit-cycle solution has the distribution of non-zero amplitudes on higher-order harmonics of $\hat{\psi}_{11}$ (e.g. $\hat{\psi}_{13}$, $\hat{\psi}_{15}$, $\cdots$, $\hat{\psi}_{31}$, $\hat{\psi}_{51}$, $\cdots$).
On the other hand, the DNS solution with the limit torus trajectory as shown in Fig.~\ref{fig:cases_XYZ}(a) does not have a similar distribution of $\hat{\psi}$ as displayed in Fig.~\ref{fig:cases_L10M10YZ}(c) but the amplitudes of other higher-order harmonics are also amplified (not shown in this paper but is qualitatively similar to Fig.~\ref{fig:cases_L10M10YZ}d).
Although the {GELE} solution considers perfect nonlinear modal interactions among the harmonics inside the domain with $l\leq 10$ and $m\leq 10$, we conjecture that {GELE} may require higher-order harmonic terms of orders $l>10$ and $m>10$ to fully reproduce the DNS solution.
We also conjecture that the DNS induces the amplification of other harmonics (e.g. $\hat{\psi}_{21}$, $\hat{\psi}_{12}$, $\cdots$) as the solutions computed in the physical space $(x,z)$ can introduce small amplitude in the non-relevant harmonics as a result of the numerical discretization.
To validate this speculation, we compute the {GELE} solution with a different initial condition where $(X,Y,Z)=(0.01,0,r-1)$ and other variables are now non-zero and random with very small initial amplitudes of order $|\hat{\psi}_{lm}|<10^{-4}$.
We clearly see in Fig.~\ref{fig:cases_L10M10YZ}(b) that the {GELE} solution with the random initial condition now exhibits a limit torus behavior after the transient period. 
It is also verified in Fig.~\ref{fig:cases_L10M10YZ}(d) that every harmonics of the {GELE} solution on the limit torus is now amplified and this amplitude distribution $\hat{\psi}_{lm}$ of the {GELE} solution resembles qualitatively the distribution of the DNS solution.


\subsection{Initial condition dependency}
\begin{figure}
\includegraphics[height=5.3cm]{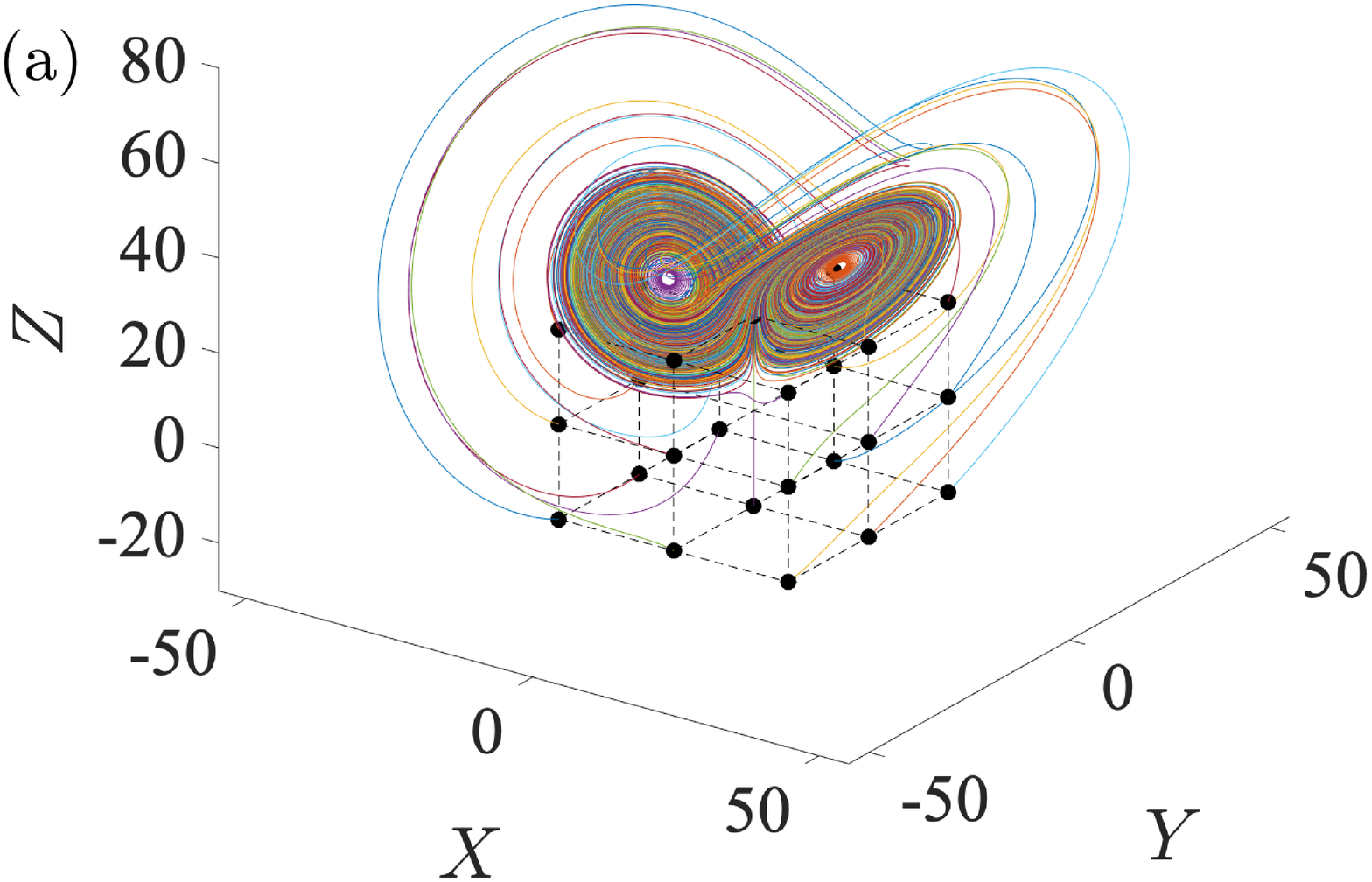}
\includegraphics[height=5.3cm]{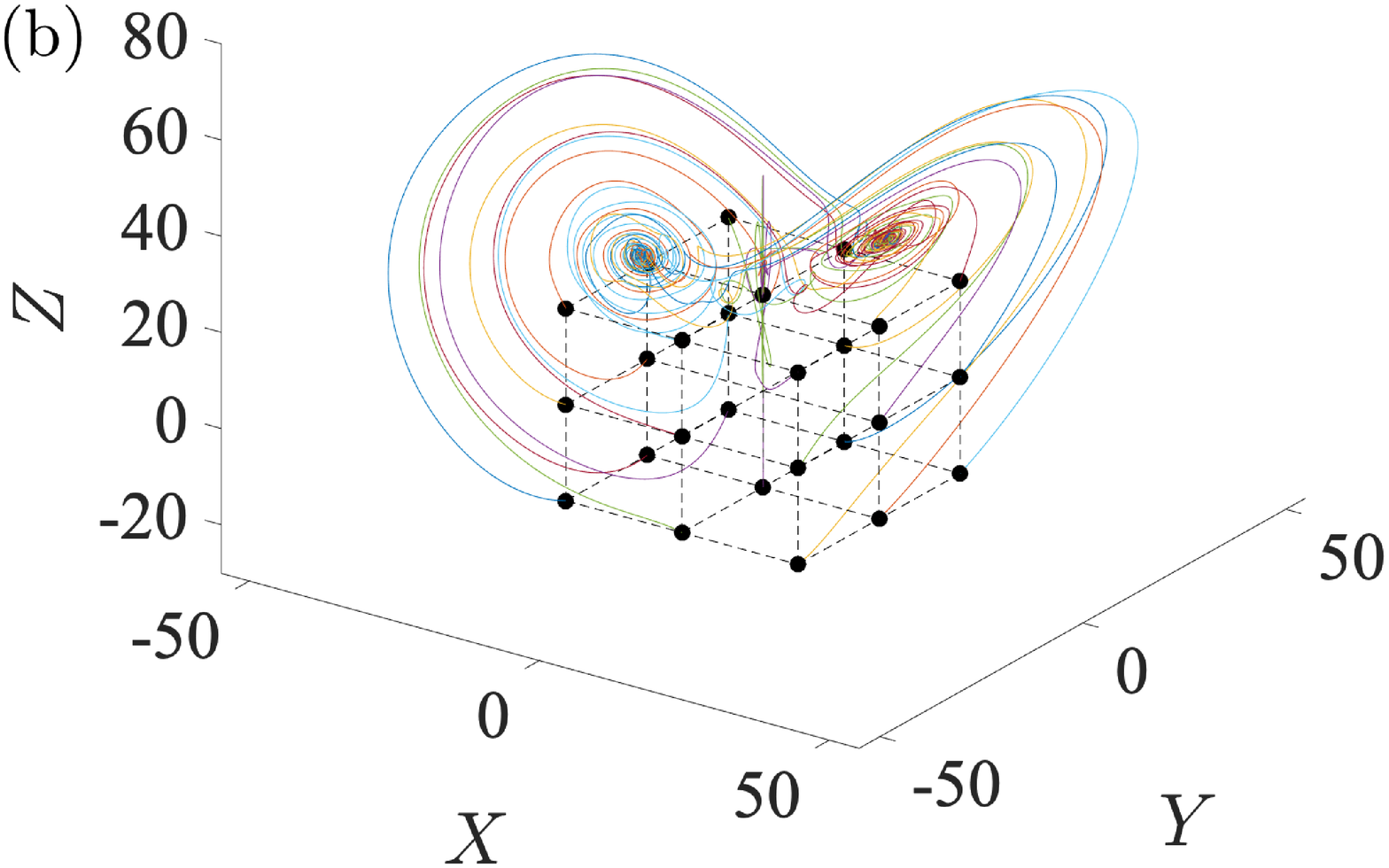}
\caption{\label{fig:IC_Lorenz} 
Trajectories on the $(X,Y,Z)$-space for the (a) Lorenz and (b) DNS solutions at $r=30$ (color solid lines). 
Black dots denote different initial conditions and dashed lines are drawn for the purpose of clear display of the initial conditions.
}
\end{figure} 
\begin{figure}
\includegraphics[height=6cm]{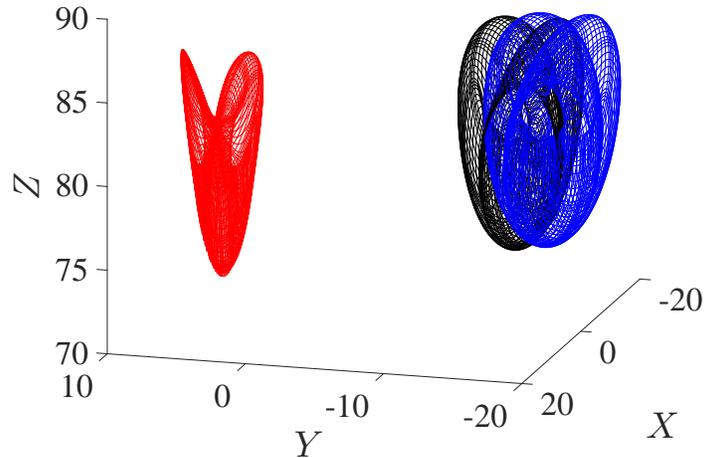}
\caption{\label{fig:IC_random} 
Trajectories on the $(X,Y,Z)$-space for DNS solutions at $r=80$ with different initial random perturbations $|\psi(x,z)|<\epsilon$ and $|\theta(x,z)|<\epsilon$ where $\epsilon=10^{-6}$ (black), $\epsilon=10^{-4}$ (blue), and $\epsilon=10^{-2}$ (red). 
}
\end{figure} 
It is now clear that the solution behavior strongly depends on the mode limits $(L,M)$ of the system, and the Lorenz equations is far different from the DNS in terms of the bifurcation behavior in the parameter space along $r$.
Other than the control parameters $(L,M)$, the initial condition also affects the bifurcation behavior since high-order systems possess multiple stable/unstable fixed points and the system's limiting dynamics can depend on the initial condition. 
As an example, we try different Lorenz-like initial conditions for the DNS and Lorenz solutions in Fig.~\ref{fig:IC_Lorenz}.
Black dots denote 26 different initial conditions generated through combinations of possible initial values $X\in\{-20, 0, 20\}$, $Y\in\{-20, 0, 20\}$ and $Z\in\{-20,0,20\}$ excluding the zero initial condition $X=Y=Z=0$. 
We see in Fig.~\ref{fig:IC_Lorenz}(a) that the Lorenz solutions at $r=30$ are chaotic and they all lie on a chaotic attractor after some transient periods. 
On the other hand, each DNS solution at $r=30$ reaches an equilibrium state and different initial conditions lead to different fixed points. 

At higher $r$, the initial condition dependency becomes more complex. 
For instance, in Fig.~\ref{fig:IC_random}, we show the DNS solutions at $r=80$ computed from initial random perturbations that satisfy $|\psi(x,z)|<\epsilon$ and $|\theta(x,z)|<\epsilon$ where $\epsilon$ is the amplitude. 
It is found that the limit tori have similar shapes for all DNS solutions, but their locations in the $(X,Y,Z)$-space vary depending on the initial amplitude $\epsilon$. 
{One difference from the Lorenz equations is that, while the Lorenz system has three fixed points, $(X,Y,Z)=(0,0,0)$ and $(X,Y,Z)=(\pm\sqrt{b(r-1)},\pm\sqrt{b(r-1)},r-1)$, the higher-order dynamical systems or the full 2D Rayleigh-B\'enard system can have many more or infinitely many fixed points, making them difficult to locate analytically.
As a result of having many fixed points, limit tori from different DNS solutions are centered at various different locations depending on the initial amplitude of perturbation. 
This is different from the Lorenz attractor, which move around the two locally unstable fixed points $(X,Y,Z)=(\pm\sqrt{b(r-1)},\pm\sqrt{b(r-1)},r-1)$. In this paper, we stop short of a full-fledged investigation of the initial condition dependency problem. It is possible, however, that the DNS solutions may possess additional fixed points with different characteristics leading to interesting conclusions; as such, the problem of multistability in DNS solutions deserves further attention in a future study.}

\section{Conclusion and discussion}
\label{sec:Conclusion}
In this paper, we propose the {generalized expansion of the Lorenz equations (GELE)} for the two-dimensional convection system, which is a generalized version of the Lorenz equations by considering higher-order harmonics in both the horizontal and vertical directions. 
{GELE} allows us to study how solutions transition from the Lorenz equations to the two-dimensional Direct Numerical Simulation (DNS) as the system orders $L$ and $M$ in the horizontal and vertical directions are varied. 
We also derived mathematical formulations for a direct comparison between the Lorenz equations, {GELE}, and DNS, and we verified in both qualitative and quantitative aspects how the Lorenz solutions in the chaotic regime are different from the DNS and high-order {GELE} solutions, which reach different equilibrium or chaotic states. 
More specifically, it is shown how the {GELE} solutions vary with $(L,M)$ and converge to those of the DNS when $L$ and $M$ are sufficiently large. 
In this study, nonlinear interactions among high-order harmonics as well as energy relations of the solutions are thoroughly analyzed. 
Furthermore, the parametric study demonstrates how trajectories of the DNS and {GELE} solutions converge to fixed points, lies on limit cycles or limit tori, depart from regular limit solutions and eventually becomes chaotic as $r$ increases. 
The initial-condition dependency is also checked to see how the {GELE} and DNS solutions behave with different initial conditions. 

The classic Lorenz equations have been considered as the minimal model that represents the chaotic nature of convection systems or even a bigger and more complex systems such as weather. 
In this study, we loosen an assumption on the minimal model by considering higher-order harmonics. 
We show by simples measures of mode amplitudes that such added complexities can lead to very different dynamical behaviors. 
The current work analyzes differences and similarities between the Lorenz equations and high-order {GELE} in a direct manner. 
And this kind of analysis should be further extended to the three-dimensional convection system to see how the increase in the spatial dimension will modify behaviors of bifurcation and chaos as the Rayleigh number increases, which will be of great interest in relevant scientific disciplines.

\section*{Supplementary Material}
{In the Supplementary Material, we demonstrate a direct comparison between the DNS and Lorenz equations by displaying the time-varying solutions of $\psi$, $\psi^{(\mathrm{Lo})}$, $\theta$, and $\theta^{(\mathrm{Lo})}$ on the plane $(x,z)$ over one streamwise wavelength $l_{x}$ for $r=30$ and $r=80$. 
In the movie, the variables $X$ and $Z$ for the DNS and Lorenz solutions are also compared. 
For $r=30$, it is clearly seen that the DNS solution reaches the equilibrium after $t>0.5$ while the Lorenz solution demonstrates a chaotic behavior. 
The chaotic variation of $X(t)$ of the Lorenz solution results in alternating appearances of positive and negative $\psi$, while the chaotic variations of $Z(t)$ and $Y(t)$ (not shown) of the Lorenz solution lead to a meandering motion in the lateral $x$-direction of $\theta$. 
It is also notable that both solutions resemble at the early development stage, but then the DNS solution deviates from the Lorenz solution as it involves nonlinear interactions among higher-order modes and reaches the steady-state equilibrium as $t$ increases.}

{For $r=80$, the DNS results of $\psi$ and $\theta$ show a more complex time-varying behavior than those at $r=30$.
For instance, at an early stage in the range $0<t<1.5$, we see a swirling motion of $\psi$ and time-periodic convective motion of $\theta$.
In the range $1.5<t<2.3$, the periodic convective motion of $\theta$ changes as the swirling motion of $\psi$ is modified in a way that the peaks of $\psi$ rotate in a wider area of the plane $(x,z)$. 
For $t>2.3$, the convective motion of $\theta$ involves lateral meandering motion and the shapes of positive/negative patches of $\psi$ become irregular. 
We note that the limit torus in Fig.~\ref{fig:cases_XYZ}(a) appears for $t>2.3$ thus we conjecture that the complex irregular motions of $\psi$ and $\theta$ with multiple time-periodicities appear as the limit torus in the phase space $(X,Y,Z)$.
The Lorenz solution at $r=80$ demonstrates a chaotic behavior in a similar manner as the Lorenz solution at $r=30$. 
}

\begin{acknowledgments}
This work was partially supported by the Small Grant for Exploratory Research (SGER) program under the National Research Foundation of Korea (NRF-2018R1D1A1A02086007).
\end{acknowledgments}

\section*{Data Availablity}
The data that support the findings of this study are available from the corresponding author upon request.

\appendix

\section{Details on convolution terms}
\label{sec:Appendix_nonlinear}
The nonlinear terms in the primitive equations (\ref{eq:psi}) and (\ref{eq:theta}):
\begin{equation}
\label{eq:Npsi}
N^{\psi}=\frac{\partial\psi}{\partial z}\frac{\partial\nabla^{2}\psi}{\partial x}-\frac{\partial\psi}{\partial x}\frac{\partial \nabla^{2}\psi}{\partial z}=\sum_{l=-L}^{L}\tilde{N}^{\psi}_{l}\exp(\mathrm{i}\alpha_{l}x),
\end{equation}
\begin{equation}
\label{eq:Ntheta}
N^{\theta}=\frac{\partial\psi}{\partial z}\frac{\partial\theta}{\partial x}-\frac{\partial\psi}{\partial x}\frac{\partial \theta}{\partial z}=\sum_{l=-L}^{L}\tilde{N}^{\theta}_{l}\exp(\mathrm{i}\alpha_{l}x),
\end{equation}
can be transformed into $\tilde{N}^{\psi}_{l}$ and $\tilde{N}^{\theta}_{l}$ that satisfy the relation (\ref{eq:nonlinear_dns}).
These nonlinear terms can be further expanded when we consider
\begin{equation}
\label{eq:tilde_N}
\tilde{N}^{\psi}_{l}=\sum_{m=0}^{M}\hat{N}^{\psi}_{lm}\sin(\beta_{m}z),~
\tilde{N}^{\theta}_{l}=\sum_{m=0}^{M}\hat{N}^{\theta}_{lm}\sin(\beta_{m}z).
\end{equation}
In the convolution process for the sine function series, we consider the relation
\begin{eqnarray}
&&\sum_{n=0}^{M}a_{n}\sin(\beta_{n}z)\sum_{k=0}^{M}b_{k}\cos(\beta_{k}z)\nonumber\\
&&=\sum_{m=0}^{M}\sum_{k=0}^{M}\left(\frac{a_{m-k}-a_{k-m}+a_{m+k}}{2}\right)b_{k}\sin(\beta_{m}z),
\end{eqnarray}
which is satisfied when we consider $a_{i}=b_{i}=0$ for indices $i<0$ or $i>M$.
Then, we get the following relations for $\hat{N}^{\psi}_{lm}$ and $\hat{N}^{\theta}_{lm}$:
\begin{eqnarray}
\label{eq:hat_Npsi}
\hat{N}^{\psi}_{lm}
&=&\sum_{j=-L}^{L}\sum_{k=0}^{M}\frac{\mathrm{i}\alpha_{j}\beta_{k}}{2}\left[\left(\alpha_{l-j}^{2}-\alpha_{j}^{2}+\beta_{k}^{2}-\beta_{m-k}^{2}\right)\hat{\psi}_{j(m-k)}\right.\nonumber\\
&&-\left(\alpha_{l-j}^{2}-\alpha_{j}^{2}+\beta_{k}^{2}-\beta_{k-m}^{2}\right)\hat{\psi}_{j(k-m)}\nonumber\\
&&\left.+\left(\alpha_{l-j}^{2}-\alpha_{j}^{2}+\beta_{k}^{2}-\beta_{m+k}^{2}\right)\hat{\psi}_{j(m+k)}\right]\hat{\psi}_{(l-j)k},
\end{eqnarray}
\begin{eqnarray}
\label{eq:hat_Ntheta}
\hat{N}^{\theta}_{l}&=&\sum_{j=-L}^{L}\sum_{k=0}^{M}\frac{\mathrm{i}\alpha_{j}\beta_{k}}{2}\left[\left(\hat{\theta}_{j(m-k)}-\hat{\theta}_{j(k-m)}+\hat{\theta}_{j(m+k)}\right)\hat{\psi}_{(l-j)k}\right.\nonumber\\
&&\left.-\left(\hat{\psi}_{j(m-k)}-\hat{\psi}_{j(k-m)}+\hat{\psi}_{j(m+k)}\right)\hat{\theta}_{(l-j)k}\right].
\end{eqnarray}

\bibliography{aipsamp}

\providecommand{\noopsort}[1]{}\providecommand{\singleletter}[1]{#1}%
\begin{thebibliography}{24}%
\makeatletter
\providecommand \@ifxundefined [1]{%
 \@ifx{#1\undefined}
}%
\providecommand \@ifnum [1]{%
 \ifnum #1\expandafter \@firstoftwo
 \else \expandafter \@secondoftwo
 \fi
}%
\providecommand \@ifx [1]{%
 \ifx #1\expandafter \@firstoftwo
 \else \expandafter \@secondoftwo
 \fi
}%
\providecommand \natexlab [1]{#1}%
\providecommand \enquote  [1]{``#1''}%
\providecommand \bibnamefont  [1]{#1}%
\providecommand \bibfnamefont [1]{#1}%
\providecommand \citenamefont [1]{#1}%
\providecommand \href@noop [0]{\@secondoftwo}%
\providecommand \href [0]{\begingroup \@sanitize@url \@href}%
\providecommand \@href[1]{\@@startlink{#1}\@@href}%
\providecommand \@@href[1]{\endgroup#1\@@endlink}%
\providecommand \@sanitize@url [0]{\catcode `\\12\catcode `\$12\catcode
  `\&12\catcode `\#12\catcode `\^12\catcode `\_12\catcode `\%12\relax}%
\providecommand \@@startlink[1]{}%
\providecommand \@@endlink[0]{}%
\providecommand \url  [0]{\begingroup\@sanitize@url \@url }%
\providecommand \@url [1]{\endgroup\@href {#1}{\urlprefix }}%
\providecommand \urlprefix  [0]{URL }%
\providecommand \Eprint [0]{\href }%
\providecommand \doibase [0]{http://dx.doi.org/}%
\providecommand \selectlanguage [0]{\@gobble}%
\providecommand \bibinfo  [0]{\@secondoftwo}%
\providecommand \bibfield  [0]{\@secondoftwo}%
\providecommand \translation [1]{[#1]}%
\providecommand \BibitemOpen [0]{}%
\providecommand \bibitemStop [0]{}%
\providecommand \bibitemNoStop [0]{.\EOS\space}%
\providecommand \EOS [0]{\spacefactor3000\relax}%
\providecommand \BibitemShut  [1]{\csname bibitem#1\endcsname}%
\let\auto@bib@innerbib\@empty
\bibitem [{\citenamefont {Getling}(1998)}]{Getling1998}%
  \BibitemOpen
  \bibfield  {author} {\bibinfo {author} {\bibfnamefont {A.~V.}\ \bibnamefont
  {Getling}},\ }\href@noop {} {\emph {\bibinfo {title} {{Rayleigh-B\'enard
  Convection: Structures and Dynamics}}}}\ (\bibinfo  {publisher} {World
  Scientific},\ \bibinfo {year} {1998})\BibitemShut {NoStop}%
\bibitem [{\citenamefont {Bodenschatz}, \citenamefont {Pesch},\ and\
  \citenamefont {Ahlers}(2000)}]{Bodenschatz2000}%
  \BibitemOpen
  \bibfield  {author} {\bibinfo {author} {\bibfnamefont {E.}~\bibnamefont
  {Bodenschatz}}, \bibinfo {author} {\bibfnamefont {W.}~\bibnamefont {Pesch}},
  \ and\ \bibinfo {author} {\bibfnamefont {G.}~\bibnamefont {Ahlers}},\
  }\bibfield  {title} {\enquote {\bibinfo {title} {{Recent developments in
  Rayleigh-B\'enard convection}},}\ }\href@noop {} {\bibfield  {journal}
  {\bibinfo  {journal} {Annu. Rev. Fluid Mech.}\ }\textbf {\bibinfo {volume}
  {32}},\ \bibinfo {pages} {709--778} (\bibinfo {year} {2000})}\BibitemShut
  {NoStop}%
\bibitem [{\citenamefont {Saltzman}(1962)}]{Saltzman1962}%
  \BibitemOpen
  \bibfield  {author} {\bibinfo {author} {\bibfnamefont {B.}~\bibnamefont
  {Saltzman}},\ }\bibfield  {title} {\enquote {\bibinfo {title} {{Finite
  amplitude free convection as an initial value problem---I}},}\ }\href@noop {}
  {\bibfield  {journal} {\bibinfo  {journal} {J. Atmos. Sci.}\ }\textbf
  {\bibinfo {volume} {19}},\ \bibinfo {pages} {329--341} (\bibinfo {year}
  {1962})}\BibitemShut {NoStop}%
\bibitem [{\citenamefont {Lorenz}(1963)}]{Lorenz1963}%
  \BibitemOpen
  \bibfield  {author} {\bibinfo {author} {\bibfnamefont {E.~N.}\ \bibnamefont
  {Lorenz}},\ }\bibfield  {title} {\enquote {\bibinfo {title} {Deterministic
  nonperiodic flow},}\ }\href@noop {} {\bibfield  {journal} {\bibinfo
  {journal} {J. Atmos. Sci.}\ }\textbf {\bibinfo {volume} {20}},\ \bibinfo
  {pages} {130--141} (\bibinfo {year} {1963})}\BibitemShut {NoStop}%
\bibitem [{\citenamefont {Tucker}(1999)}]{tucker1999lorenz}%
  \BibitemOpen
  \bibfield  {author} {\bibinfo {author} {\bibfnamefont {W.}~\bibnamefont
  {Tucker}},\ }\bibfield  {title} {\enquote {\bibinfo {title} {The {L}orenz
  attractor exists},}\ }\href@noop {} {\bibfield  {journal} {\bibinfo
  {journal} {C. R. Acad. Sci.---S\'{e}r. {I}---Math.}\ }\textbf {\bibinfo
  {volume} {328}},\ \bibinfo {pages} {1197--1202} (\bibinfo {year}
  {1999})}\BibitemShut {NoStop}%
\bibitem [{\citenamefont {Gleick}(1987)}]{gleick1988chaos}%
  \BibitemOpen
  \bibfield  {author} {\bibinfo {author} {\bibfnamefont {J.}~\bibnamefont
  {Gleick}},\ }\href@noop {} {\emph {\bibinfo {title} {Chaos: Making a New
  Science}}}\ (\bibinfo  {publisher} {Viking Penguin},\ \bibinfo {address} {New
  York},\ \bibinfo {year} {1987})\ p.\ \bibinfo {pages} {400}\BibitemShut
  {NoStop}%
\bibitem [{\citenamefont {Shen}(2014)}]{Shen2014}%
  \BibitemOpen
  \bibfield  {author} {\bibinfo {author} {\bibfnamefont {B.~W.}\ \bibnamefont
  {Shen}},\ }\bibfield  {title} {\enquote {\bibinfo {title} {{{Nonlinear
  feedback in a five-dimensional Lorenz model}}},}\ }\href@noop {} {\bibfield
  {journal} {\bibinfo  {journal} {J. Atmos. Sci.}\ }\textbf {\bibinfo {volume}
  {71}},\ \bibinfo {pages} {1701--1723} (\bibinfo {year} {2014})}\BibitemShut
  {NoStop}%
\bibitem [{\citenamefont {Stenflo}(1996)}]{Stenflo1996}%
  \BibitemOpen
  \bibfield  {author} {\bibinfo {author} {\bibfnamefont {L.}~\bibnamefont
  {Stenflo}},\ }\bibfield  {title} {\enquote {\bibinfo {title} {Generalized
  {L}orenz equations for acoustic-gravity waves in the atmosphere},}\
  }\href@noop {} {\bibfield  {journal} {\bibinfo  {journal} {Phys. Scr.}\
  }\textbf {\bibinfo {volume} {53}},\ \bibinfo {pages} {83--84} (\bibinfo
  {year} {1996})}\BibitemShut {NoStop}%
\bibitem [{\citenamefont {Park}\ \emph
  {et~al.}(2015{\natexlab{a}})\citenamefont {Park}, \citenamefont {Han},
  \citenamefont {Lee}, \citenamefont {Jeon},\ and\ \citenamefont
  {Baik}}]{Park2016PS}%
  \BibitemOpen
  \bibfield  {author} {\bibinfo {author} {\bibfnamefont {J.}~\bibnamefont
  {Park}}, \bibinfo {author} {\bibfnamefont {B.-S.}\ \bibnamefont {Han}},
  \bibinfo {author} {\bibfnamefont {H.}~\bibnamefont {Lee}}, \bibinfo {author}
  {\bibfnamefont {Y.-L.}\ \bibnamefont {Jeon}}, \ and\ \bibinfo {author}
  {\bibfnamefont {J.-J.}\ \bibnamefont {Baik}},\ }\bibfield  {title} {\enquote
  {\bibinfo {title} {{Stability and periodicity of high-order Lorenz-Stenflo
  equations}},}\ }\href@noop {} {\bibfield  {journal} {\bibinfo  {journal}
  {Phys. Scr.}\ }\textbf {\bibinfo {volume} {91}},\ \bibinfo {pages} {065202}
  (\bibinfo {year} {2015}{\natexlab{a}})}\BibitemShut {NoStop}%
\bibitem [{\citenamefont {Moon}\ \emph {et~al.}(2019)\citenamefont {Moon},
  \citenamefont {Seo}, \citenamefont {Han}, \citenamefont {Park},\ and\
  \citenamefont {Baik}}]{Moon2019}%
  \BibitemOpen
  \bibfield  {author} {\bibinfo {author} {\bibfnamefont {S.}~\bibnamefont
  {Moon}}, \bibinfo {author} {\bibfnamefont {J.~M.}\ \bibnamefont {Seo}},
  \bibinfo {author} {\bibfnamefont {B.-S.}\ \bibnamefont {Han}}, \bibinfo
  {author} {\bibfnamefont {J.}~\bibnamefont {Park}}, \ and\ \bibinfo {author}
  {\bibfnamefont {J.-J.}\ \bibnamefont {Baik}},\ }\bibfield  {title} {\enquote
  {\bibinfo {title} {{{A physically extended Lorenz system}}},}\ }\href@noop {}
  {\bibfield  {journal} {\bibinfo  {journal} {Chaos}\ }\textbf {\bibinfo
  {volume} {29}},\ \bibinfo {pages} {063129} (\bibinfo {year}
  {2019})}\BibitemShut {NoStop}%
\bibitem [{\citenamefont {Felicio}\ and\ \citenamefont
  {Rech}(2018)}]{Felicio2018}%
  \BibitemOpen
  \bibfield  {author} {\bibinfo {author} {\bibfnamefont {C.~C.}\ \bibnamefont
  {Felicio}}\ and\ \bibinfo {author} {\bibfnamefont {P.~C.}\ \bibnamefont
  {Rech}},\ }\bibfield  {title} {\enquote {\bibinfo {title} {{{On the dynamics
  of five- and six-dimensional Lorenz models}}},}\ }\href@noop {} {\bibfield
  {journal} {\bibinfo  {journal} {J. Phys. Commun.}\ }\textbf {\bibinfo
  {volume} {2}},\ \bibinfo {pages} {025028} (\bibinfo {year}
  {2018})}\BibitemShut {NoStop}%
\bibitem [{\citenamefont {Moon}\ \emph {et~al.}(2017)\citenamefont {Moon},
  \citenamefont {Han}, \citenamefont {Park}, \citenamefont {Seo},\ and\
  \citenamefont {Baik}}]{Moon2017}%
  \BibitemOpen
  \bibfield  {author} {\bibinfo {author} {\bibfnamefont {S.}~\bibnamefont
  {Moon}}, \bibinfo {author} {\bibfnamefont {B.-S.}\ \bibnamefont {Han}},
  \bibinfo {author} {\bibfnamefont {J.}~\bibnamefont {Park}}, \bibinfo {author}
  {\bibfnamefont {J.~M.}\ \bibnamefont {Seo}}, \ and\ \bibinfo {author}
  {\bibfnamefont {J.-J.}\ \bibnamefont {Baik}},\ }\bibfield  {title} {\enquote
  {\bibinfo {title} {{{Periodicity and chaos of high-order Lorenz systems}}},}\
  }\href@noop {} {\bibfield  {journal} {\bibinfo  {journal} {Int. J.
  Bifurcation Chaos}\ }\textbf {\bibinfo {volume} {27}},\ \bibinfo {pages}
  {1750176} (\bibinfo {year} {2017})}\BibitemShut {NoStop}%
\bibitem [{\citenamefont {Moon}, \citenamefont {Seo},\ and\ \citenamefont
  {Baik}(2020)}]{Moon2020}%
  \BibitemOpen
  \bibfield  {author} {\bibinfo {author} {\bibfnamefont {S.}~\bibnamefont
  {Moon}}, \bibinfo {author} {\bibfnamefont {J.~M.}\ \bibnamefont {Seo}}, \
  and\ \bibinfo {author} {\bibfnamefont {J.-J.}\ \bibnamefont {Baik}},\
  }\bibfield  {title} {\enquote {\bibinfo {title} {{{High-dimensional
  generalizations of the Lorenz system and implications for
  predictability}}},}\ }\href@noop {} {\bibfield  {journal} {\bibinfo
  {journal} {Phys. Scr.}\ }\textbf {\bibinfo {volume} {95}},\ \bibinfo {pages}
  {085209} (\bibinfo {year} {2020})}\BibitemShut {NoStop}%
\bibitem [{\citenamefont {Stevens}(2011)}]{Stevens2011}%
  \BibitemOpen
  \bibfield  {author} {\bibinfo {author} {\bibfnamefont {R.}~\bibnamefont
  {Stevens}},\ }\emph {\bibinfo {title} {{Rayleigh-B\'enard Turbulence}}},\
  \href@noop {} {Ph.D. thesis},\ \bibinfo  {school} {University of Twente}
  (\bibinfo {year} {2011})\BibitemShut {NoStop}%
\bibitem [{\citenamefont {Bao}, \citenamefont {Luo},\ and\ \citenamefont
  {Ye}(2017)}]{Bao2017}%
  \BibitemOpen
  \bibfield  {author} {\bibinfo {author} {\bibfnamefont {Y.}~\bibnamefont
  {Bao}}, \bibinfo {author} {\bibfnamefont {J.}~\bibnamefont {Luo}}, \ and\
  \bibinfo {author} {\bibfnamefont {M.}~\bibnamefont {Ye}},\ }\bibfield
  {title} {\enquote {\bibinfo {title} {{{Parallel direct method of DNS for
  two-dimensional turbulent Rayleigh-B\'enard convection}}},}\ }\href@noop {}
  {\bibfield  {journal} {\bibinfo  {journal} {J. Mech.}\ }\textbf {\bibinfo
  {volume} {34}},\ \bibinfo {pages} {159--166} (\bibinfo {year}
  {2017})}\BibitemShut {NoStop}%
\bibitem [{\citenamefont {Paul}\ \emph {et~al.}(2012)\citenamefont {Paul},
  \citenamefont {Verma}, \citenamefont {Wahi}, \citenamefont {Reddy},\ and\
  \citenamefont {Kumar}}]{Paul2012}%
  \BibitemOpen
  \bibfield  {author} {\bibinfo {author} {\bibfnamefont {S.}~\bibnamefont
  {Paul}}, \bibinfo {author} {\bibfnamefont {M.~K.}\ \bibnamefont {Verma}},
  \bibinfo {author} {\bibfnamefont {P.}~\bibnamefont {Wahi}}, \bibinfo {author}
  {\bibfnamefont {S.~K.}\ \bibnamefont {Reddy}}, \ and\ \bibinfo {author}
  {\bibfnamefont {K.}~\bibnamefont {Kumar}},\ }\bibfield  {title} {\enquote
  {\bibinfo {title} {{{Bifurcation analysis of the flow patterns in
  two-dimensional Rayleigh-B\'enard convection}}},}\ }\href@noop {} {\bibfield
  {journal} {\bibinfo  {journal} {Int. J. Bifurcation Chaos}\ }\textbf
  {\bibinfo {volume} {22}},\ \bibinfo {pages} {1230018} (\bibinfo {year}
  {2012})}\BibitemShut {NoStop}%
\bibitem [{\citenamefont {Weideman}\ and\ \citenamefont
  {Reddy}(2000)}]{Weideman2000}%
  \BibitemOpen
  \bibfield  {author} {\bibinfo {author} {\bibfnamefont {J.~A.}\ \bibnamefont
  {Weideman}}\ and\ \bibinfo {author} {\bibfnamefont {S.~C.}\ \bibnamefont
  {Reddy}},\ }\bibfield  {title} {\enquote {\bibinfo {title} {{A MATLAB
  differentiation matrix suite}},}\ }\href@noop {} {\bibfield  {journal}
  {\bibinfo  {journal} {ACM Trans. Math. Softw.}\ }\textbf {\bibinfo {volume}
  {26}},\ \bibinfo {pages} {465--519} (\bibinfo {year} {2000})}\BibitemShut
  {NoStop}%
\bibitem [{\citenamefont {Antkowiak}(2005)}]{Antkowiak2005}%
  \BibitemOpen
  \bibfield  {author} {\bibinfo {author} {\bibfnamefont {A.}~\bibnamefont
  {Antkowiak}},\ }\emph {\bibinfo {title} {Dynamique aux temps courts d'un
  tourbillon isol\'e}},\ \href@noop {} {Ph.D. thesis},\ \bibinfo  {school}
  {Universit\'e Paul Sabatier de Toulouse} (\bibinfo {year} {2005})\BibitemShut
  {NoStop}%
\bibitem [{\citenamefont {Park}, \citenamefont {Billant},\ and\ \citenamefont
  {Baik}(2017)}]{Park2017}%
  \BibitemOpen
  \bibfield  {author} {\bibinfo {author} {\bibfnamefont {J.}~\bibnamefont
  {Park}}, \bibinfo {author} {\bibfnamefont {P.}~\bibnamefont {Billant}}, \
  and\ \bibinfo {author} {\bibfnamefont {J.-J.}\ \bibnamefont {Baik}},\
  }\bibfield  {title} {\enquote {\bibinfo {title} {{Instabilities and transient
  growth of the stratified Taylor-Couette flow in a Rayleigh-unstable
  regime}},}\ }\href@noop {} {\bibfield  {journal} {\bibinfo  {journal} {J.
  Fluid Mech.}\ }\textbf {\bibinfo {volume} {822}},\ \bibinfo {pages} {80--108}
  (\bibinfo {year} {2017})}\BibitemShut {NoStop}%
\bibitem [{\citenamefont {Kim}, \citenamefont {Moin},\ and\ \citenamefont
  {Moser}(1987)}]{KMM1987}%
  \BibitemOpen
  \bibfield  {author} {\bibinfo {author} {\bibfnamefont {J.}~\bibnamefont
  {Kim}}, \bibinfo {author} {\bibfnamefont {P.}~\bibnamefont {Moin}}, \ and\
  \bibinfo {author} {\bibfnamefont {R.}~\bibnamefont {Moser}},\ }\bibfield
  {title} {\enquote {\bibinfo {title} {{{Turbulence statistics in fully
  developed channel flow at low Reynolds number}}},}\ }\href@noop {} {\bibfield
   {journal} {\bibinfo  {journal} {J. Fluid Mech.}\ }\textbf {\bibinfo {volume}
  {177}},\ \bibinfo {pages} {133--166} (\bibinfo {year} {1987})}\BibitemShut
  {NoStop}%
\bibitem [{\citenamefont {Yu}, \citenamefont {Zhou},\ and\ \citenamefont
  {Lai}(1996)}]{Yu1996}%
  \BibitemOpen
  \bibfield  {author} {\bibinfo {author} {\bibfnamefont {M.~Y.}\ \bibnamefont
  {Yu}}, \bibinfo {author} {\bibfnamefont {C.~T.}\ \bibnamefont {Zhou}}, \ and\
  \bibinfo {author} {\bibfnamefont {C.~H.}\ \bibnamefont {Lai}},\ }\bibfield
  {title} {\enquote {\bibinfo {title} {{{The bifurcation characteristics of the
  generalized Lorenz equations}}},}\ }\href@noop {} {\bibfield  {journal}
  {\bibinfo  {journal} {Phys. Scr.}\ }\textbf {\bibinfo {volume} {53}},\
  \bibinfo {pages} {321} (\bibinfo {year} {1996})}\BibitemShut {NoStop}%
\bibitem [{\citenamefont {Park}\ \emph
  {et~al.}(2015{\natexlab{b}})\citenamefont {Park}, \citenamefont {Lee},
  \citenamefont {Jeon},\ and\ \citenamefont {Baik}}]{Park2015PS}%
  \BibitemOpen
  \bibfield  {author} {\bibinfo {author} {\bibfnamefont {J.}~\bibnamefont
  {Park}}, \bibinfo {author} {\bibfnamefont {H.}~\bibnamefont {Lee}}, \bibinfo
  {author} {\bibfnamefont {Y.-L.}\ \bibnamefont {Jeon}}, \ and\ \bibinfo
  {author} {\bibfnamefont {J.-J.}\ \bibnamefont {Baik}},\ }\bibfield  {title}
  {\enquote {\bibinfo {title} {{Periodicity of the Lorenz-Stenflo
  equations}},}\ }\href@noop {} {\bibfield  {journal} {\bibinfo  {journal}
  {Phys. Scr.}\ }\textbf {\bibinfo {volume} {90}},\ \bibinfo {pages} {065201}
  (\bibinfo {year} {2015}{\natexlab{b}})}\BibitemShut {NoStop}%
\bibitem [{\citenamefont {Dullin}\ \emph {et~al.}(2007)\citenamefont {Dullin},
  \citenamefont {Schmidt}, \citenamefont {Richter},\ and\ \citenamefont
  {Grossmann}}]{Dullin2007}%
  \BibitemOpen
  \bibfield  {author} {\bibinfo {author} {\bibfnamefont {H.~R.}\ \bibnamefont
  {Dullin}}, \bibinfo {author} {\bibfnamefont {S.}~\bibnamefont {Schmidt}},
  \bibinfo {author} {\bibfnamefont {P.~H.}\ \bibnamefont {Richter}}, \ and\
  \bibinfo {author} {\bibfnamefont {S.~K.}\ \bibnamefont {Grossmann}},\
  }\bibfield  {title} {\enquote {\bibinfo {title} {{{Extended phase diagram of
  the Lorenz model}}},}\ }\href@noop {} {\bibfield  {journal} {\bibinfo
  {journal} {Int. J. Bifurcation Chaos}\ }\textbf {\bibinfo {volume} {17}},\
  \bibinfo {pages} {3013--3033} (\bibinfo {year} {2007})}\BibitemShut {NoStop}%
\bibitem [{\citenamefont {Grebogi}, \citenamefont {Ott},\ and\ \citenamefont
  {Yorke}(1985)}]{Grebogi1985}%
  \BibitemOpen
  \bibfield  {author} {\bibinfo {author} {\bibfnamefont {C.}~\bibnamefont
  {Grebogi}}, \bibinfo {author} {\bibfnamefont {E.}~\bibnamefont {Ott}}, \ and\
  \bibinfo {author} {\bibfnamefont {J.~A.}\ \bibnamefont {Yorke}},\ }\bibfield
  {title} {\enquote {\bibinfo {title} {{Attractors on an $N$-torus:
  Quasiperiodicity versus chaos}},}\ }\href@noop {} {\bibfield  {journal}
  {\bibinfo  {journal} {Phys. D}\ }\textbf {\bibinfo {volume} {15}},\ \bibinfo
  {pages} {354--373} (\bibinfo {year} {1985})}\BibitemShut {NoStop}%
\end{thebibliography}%

\end{document}